\pdfoutput=1
\documentclass[10pt]{amsart}

\usepackage[bookmarks,colorlinks=false,backref=false,plainpages=false]{hyperref}
\usepackage{geometry}
\usepackage{amsmath,amssymb,amsthm}
\usepackage{graphicx}
\usepackage{booktabs}

\newtheorem{thm}{Theorem}[section]
\newtheorem{cor}[thm]{Corollary}
\newtheorem{prop}[thm]{Proposition}
\theoremstyle{remark}
\newtheorem{rem}[thm]{Remark}

\usepackage{amsaddr}

\numberwithin{equation}{section}

\renewcommand\O{\Omega}
\def\G{\Gamma}
\newcommand{\bsigma}{{\boldsymbol\sigma}}
\newcommand{\btau}{{\boldsymbol\tau}}
\newcommand{\bchi}{{\boldsymbol\chi}}
\newcommand{\brho}{{\boldsymbol\rho}}
\newcommand{\bmu}{{\boldsymbol\mu}}
\newcommand{\bSigma}{{\boldsymbol\Sigma}}
\def\bb{{\bf b}}
\def\bff{{\bf f}}
\def\bg{{\bf g}}
\def\bn{{\bf n}}
\def\br{{\bf r}}
\def\bzero{{\bf 0}}
\def\cT{{\mathcal{T}}}
\def\R{\mathbb{R}}
\newcommand{\divv}{\operatorname{div}}
\newcommand{\LS}{{\mathtt{LS}}}
\newcommand{\Rcp}{R^{\rm cst}_p}
\newcommand{\Rbp}{R^{\rm bal}_p}
\newcommand{\Rcd}{R^{\rm cst}_d}
\newcommand{\Rbd}{R^{\rm bal}_d}
\newcommand{\tri}{|\!|\!|}
\newcommand{\beq}{\begin{equation}}
\newcommand{\eeq}{\end{equation}}

\title[Signed QoI gaps for goal-oriented error estimation]{Signed QoI
Gaps for Goal-Oriented Error Estimation:\\ Reliability and Sharpness from the
Same First-Order System Residuals}
\author[S. Zhang]{Shun Zhang}
\address{Department of Mathematics, City University of Hong Kong, Kowloon Tong, Hong Kong, China}
\email{shun.zhang@cityu.edu.hk}
\thanks{This work was supported in part by
Research Grants Council of the Hong Kong SAR, China, under the GRF Grant Project No. CityU 11316222, CityU 11305025}
\date{\today}

\begin{document}

\begin{abstract}
There are two common ways to estimate the error in a quantity of interest
(QoI).
The product of primal and adjoint error estimators gives a reliable but
pessimistic bound. The dual-weighted residual method gives a sharp signed
estimate from a more accurate approximation of the adjoint solution. We
write the primal and adjoint problems as first-order systems and obtain
both a reliable bound and a sharp estimate from the same residuals.
The first paper of this series (arXiv:2607.23850) uses two
consistent evaluators of a quantity of interest, one built from the
potentials and one in which the gradients are expressed through the
fluxes. We show that they differ exactly by the $A^{-1}$ inner product of
the primal and adjoint constitutive residuals, and we call this difference
the signed QoI gap. The identity holds for arbitrary primal and
adjoint approximations. In our numerical experiments the signed QoI gap
is of the size of the larger of the errors of the two evaluators, while
the product of the estimators exceeds this error by factors of up to
several hundred. We therefore use the least-squares estimators to refine the mesh
and the signed QoI gap to decide when to stop. We also introduce a third
evaluator, built from the fluxes alone, which restores the reaction term.
The resulting reaction-restored QoI gap estimates the error when reaction
dominates. With the non-intrusive least-squares flux recovery, we can
apply the method to any discretization that
produces a conforming potential, including finite difference and finite
volume schemes.
\end{abstract}

\maketitle

\section{Introduction}\label{sec:intro}

In many applications the quantity of interest (QoI) is a scalar functional of
the solution of a partial differential equation, such as a drag
coefficient, an average over a subdomain, or a flux through an interface.
We want to know how large the error of the computed value is. In the standard Galerkin
setting, with the Galerkin approximation $u_h\in X_h$, $\gamma=m(u)$ and the adjoint solution $z$ satisfying
$a(s,z)=m(s)$ for all admissible $s$,
\[
\gamma-m(u_h)=a(u-u_h,\,z-\psi_h)\qquad\forall\,\psi_h\in X_h .
\]
The dual-weighted residual method \cite{BR:01,BR:03} replaces $z$ by a
more accurate approximation of the adjoint solution, computed in an
enriched space or reconstructed from the Galerkin approximation
$z_h\in X_h$ of $z$. The resulting signed estimate is accurate in practice
but not guaranteed, and its elementwise terms serve as refinement
indicators.
The other common approach takes $\psi_h=z_h$ and bounds the two factors by
estimators of the primal and the adjoint error \cite{PO:99,OP:01},
$|\gamma-m(u_h)|\le C\,\eta_p\eta_d$. This bound is reliable and is the
basis of the convergence theory of goal-oriented adaptive methods
\cite{MS:09,BET:11,FPV:16,IP:21}, but it can be very pessimistic. In this
paper we obtain both a reliable bound and a signed estimate from the same
first-order system residuals.

The first paper of this series \cite{WZ:26} considers the general
second-order elliptic equation with diffusion coefficient $A$, convection
coefficient $\bb$ and reaction coefficient $c$, which is also our model
problem. That paper writes the primal problem and a physically
meaningful adjoint problem as first-order systems for a flux and a
potential. The norms of their residuals, the
least-squares estimators $\eta_p^{\LS}$ and $\eta_d^{\LS}$, are equivalent
to the primal and adjoint errors for arbitrary approximations
\cite{CLMM:94,CMM:97,BG:98,BG:09,Zhang:23,CZ:10b}. The same paper uses
two consistent evaluators of $\gamma$. The first, which we call the potential
evaluator, is the corrected functional of Giles and S\"uli \cite{GS:02},
$\gamma_0(w,v):=m(w)+\ell(v)-a(w,v)$, where $\ell$ is the right-hand side
of the primal problem. It satisfies $\gamma-\gamma_0(w,v)=a(u-w,z-v)$ for
all admissible $w$ and $v$, so we do not need Galerkin orthogonality. The
second, $\gamma_1$, also uses the flux approximations that a least-squares
method computes along with the potentials. It expresses the gradients in
$\gamma_0$ through these fluxes and adds two residual corrections.
Another way to combine least squares with a quantity of interest is to
build the quantity into the least-squares functional itself
\cite{CCLMOT:14}.

The product of the two estimators gives a reliable upper bound for the
errors of both evaluators, and its elementwise contributions
\cite{WZ:26} are good refinement indicators. In our experiments, however,
the product exceeds the larger of these errors by factors of up to several
hundred. If we stop when the product falls below a tolerance, we refine far
more than needed. This overestimate does not come from the estimators.
For the Poisson model problem each of them
is asymptotically exact at the least-squares solution \cite{CS:18}, so
there the whole overestimate comes from bounding the error $a(u-w,z-v)$ of
$\gamma_0(w,v)$ by the product of the two error norms. To decide when to
stop, we need to know how large the error actually is.

For this we compare two approximations of the same quantity. We compute
$\gamma$ twice from the same primal and adjoint approximations
$(\btau,w)$ and $(\brho,v)$, once as $\gamma_0(w,v)$ from the potentials
and once as $\gamma_1(w,\btau,v,\brho)$ with the gradients expressed
through the fluxes. We do not compute a second, more accurate
approximation. We show that
the difference of the two values, which we call the {\em signed QoI
gap}, is exactly the $A^{-1}$ inner product of the primal and adjoint
constitutive residuals,
\[
\gamma_1(w,\btau,v,\brho)-\gamma_0(w,v)=\big(A^{-1}\Rcp(\btau,w),\;\Rcd(\brho,v)\big),
\]
with $\Rcp(\btau,w)=\bff_2+\btau+A\nabla w$ and
$\Rcd(\brho,v)=\bg_2+\brho+A\nabla v+\bb v$, where $\bff_2$ and $\bg_2$
are the vector-valued data of the primal problem and of the quantity of
interest. The identity is algebraic and holds for arbitrary primal and
adjoint approximations. When $\gamma$ lies
between the two values, the QoI gap is between one and two
times the larger of the two errors. Otherwise the QoI gap is smaller than
the larger error, and much smaller only when the two values have nearly
the same error.

The dual-weighted residual estimate also computes
$\gamma$ twice, but it evaluates $\gamma_0$ the second time at a more
accurate adjoint approximation. It is a good estimate when the second
value is much more accurate than the first. The
signed QoI gap is a generalized dual-weighted residual estimate with
flux-based weights. We do not assume that either of its two values is
the more accurate, and we need no more accurate approximation. The fluxes
need not be more accurate than the gradients of the potentials either.
The signed QoI gap and the product bound
come from the same residuals. For the QoI gap we take the inner product
of the primal and adjoint constitutive residuals, and for the bound we
take the product of the norms of the primal and adjoint residuals.

Since $-A^{-1}\Rcp(\btau,w)=\nabla(u-w)+A^{-1}(\bsigma-\btau)$ and
$-\Rcd(\brho,v)=A\nabla(z-v)+\bb(z-v)+(\br-\brho)$, the QoI gap contains the
error $a(u-w,z-v)$ of $\gamma_0(w,v)$ without its reaction term, plus
terms in the flux errors. Where reaction dominates, the reaction term can
dominate the errors of both evaluators, and the signed QoI gap does not
contain it. For this case we introduce the reaction-restored evaluator
$\gamma_2(\btau,\brho)$, built from the fluxes alone. The reaction-restored
QoI gap $\gamma_2(\btau,\brho)-\gamma_0(w,v)$ then estimates the error of
$\gamma_0(w,v)$. We refine the mesh with the least-squares estimators and
stop when the magnitude of the signed QoI gap, or of the
reaction-restored QoI gap where reaction dominates, is below a tolerance. In our experiments the QoI gaps follow the errors
closely, even where the product bound is very pessimistic.

The construction does not depend on the discretization. We can take the
potentials from any method that produces conforming potentials, and the
fluxes from least-squares flux recovery \cite{CZ:10b,LZYZ:24} or from a
mixed method. An augmented mixed method \cite{MH:02,LZ:24} gives both, and
so does a least-squares finite element method. Our experiments include potentials from finite difference and
finite volume schemes, equilibrated mixed fluxes for Kellogg's problem, a
comparison with a dual-weighted residual computation, and computations
with reaction coefficients from $c=1$ to $c=10^8$.

Section~\ref{sec:setting} recalls the setting and the evaluators of
\cite{WZ:26}. Section~\ref{sec:gap} derives the QoI gap identity, shows
what the signed QoI gap estimates, and introduces $\gamma_2$ and the
reaction-restored QoI gap.
Section~\ref{sec:relation} relates the QoI gaps to the product bound, the
dual-weighted residual method, flux recovery and, for symmetric problems,
the primal--dual bounds of \cite{Zhang:20}. Section~\ref{sec:routes}
describes how we compute the fluxes and how we refine the meshes, and
Section~\ref{sec:num} reports the numerical experiments.

\section{Primal and adjoint first-order systems and two evaluators}\label{sec:setting}

\subsection{The primal problem, the quantity of interest and the adjoint problem}

Let $\O\subset\R^d$, $d=2,3$, be a bounded Lipschitz polyhedral domain,
whose boundary is split into relatively open, disjoint parts $\G_D$ and
$\G_N$ with $\partial\O=\overline{\G}_D\cup\overline{\G}_N$, and let $\bn$
be the outward unit normal. The coefficients are
$A\in L^\infty(\O)^{d\times d}$, symmetric and uniformly positive definite,
$\bb\in L^\infty(\O)^d$ and $c\in L^\infty(\O)$. The data of the primal
problem are $f_1\in L^2(\O)$, $\bff_2\in L^2(\O)^d$,
$\phi_p\in H^{1/2}(\G_D)$ and $\psi_p\in L^2(\G_N)$, and those of the
quantity of interest are $g_1\in L^2(\O)$, $\bg_2\in L^2(\O)^d$,
$\phi_d\in H^{1/2}(\G_D)$ and $\psi_d\in L^2(\G_N)$. We write
$(\cdot,\cdot)$ and $\|\cdot\|_0$ for the inner product and norm of
$L^2(\O)$, and $\langle\cdot,\cdot\rangle_{\G}$ for the $L^2(\G)$ inner
product on $\G\subset\partial\O$, or the duality pairing of $H^{-1/2}$ and
$H^{1/2}$ there. The spaces are
\begin{align*}
H^1_D(\O)&:=\{s\in H^1(\O):\ s=0 \text{ on }\G_D\}, &
H_N(\divv;\O)&:=\{\bmu\in H(\divv;\O):\ \bmu\cdot\bn=0 \text{ on }\G_N\},\\
X_p&:=\{w\in H^1(\O):\ w=\phi_p \text{ on }\G_D\}, &
\bSigma_p&:=\{\btau\in H(\divv;\O):\ \btau\cdot\bn=\psi_p \text{ on }\G_N\},\\
X_d&:=\{v\in H^1(\O):\ v=\phi_d \text{ on }\G_D\}, &
\bSigma_d&:=\{\brho\in H(\divv;\O):\ \brho\cdot\bn=\psi_d \text{ on }\G_N\}.
\end{align*}

The primal problem is
\beq \label{primal_pde}
-\divv(A\nabla u)+\bb\cdot\nabla u+cu = f_1+\divv\bff_2 \ \text{ in }\O,
\qquad u=\phi_p \ \text{ on }\G_D, \qquad
-(A\nabla u+\bff_2)\cdot\bn=\psi_p \ \text{ on }\G_N .
\eeq
With
\beq \label{forms}
a(w,v):=(A\nabla w,\nabla v)+(\bb\cdot\nabla w+cw,v), \qquad
\ell(v):=(f_1,v)-(\bff_2,\nabla v)-\langle\psi_p,v\rangle_{\G_N},
\eeq
its weak form is: find $u\in X_p$ with $a(u,s)=\ell(s)$ for all
$s\in H^1_D(\O)$. We assume that $a$ satisfies the inf-sup condition on
$H^1_D(\O)\times H^1_D(\O)$, so that this problem and its adjoint problem are uniquely solvable. The flux
$\bsigma:=-(A\nabla u+\bff_2)$ belongs to $\bSigma_p$, and $(\bsigma,u)$
solves the first-order system
\beq \label{fos_p}
\divv\bsigma+\bb\cdot\nabla u+cu=f_1, \qquad
\bsigma+A\nabla u = -\bff_2 \qquad \text{in }\O,
\eeq
whose first equation we call the balance equation and the second the
constitutive equation.

The quantity of interest is
\beq \label{QoI}
\gamma:=m(u)-\langle(A\nabla u+\bff_2)\cdot\bn,\phi_d\rangle_{\G_D}
=m(u)+\langle\bsigma\cdot\bn,\phi_d\rangle_{\G_D},
\qquad
m(w):=(g_1,w)-(\bg_2,\nabla w)-\langle\psi_d,w\rangle_{\G_N}.
\eeq
The four data of the quantity of interest determine what we measure. The
datum $g_1$ gives a weighted integral of $u$ over $\O$, $\bg_2$ a weighted
integral of $\nabla u$, $\psi_d$ a weighted integral of $u$ on $\G_N$, and
$\phi_d$ a weighted integral of the normal flux $\bsigma\cdot\bn$ on
$\G_D$. For instance, with the indicator function $\mathbf 1_\omega$ of a
subdomain $\omega$, the choice $g_1=\mathbf 1_\omega$ gives the integral of
$u$ over $\omega$, and in two dimensions $\bg_2=(\mathbf 1_\omega,0)^T$ gives
$-\int_\omega\partial_xu\,dx$. The choice $\psi_d=-1$ gives the integral
of $u$ over $\G_N$, and $\phi_d=1$ gives the total outward flux
$\int_{\G_D}\bsigma\cdot\bn\,ds$ through $\G_D$.

The adjoint problem is: find $z\in X_d$ with
\beq \label{dual_weak}
a(s,z)=m(s)\qquad\forall\,s\in H^1_D(\O).
\eeq
Its strong form is the physically meaningful adjoint problem of
\cite{WZ:26},
\beq \label{dual_pde}
-\divv(A\nabla z+\bb z)+cz=g_1+\divv\bg_2 \ \text{ in }\O, \qquad
z=\phi_d \ \text{ on }\G_D, \qquad
-(A\nabla z+\bb z+\bg_2)\cdot\bn=\psi_d \ \text{ on }\G_N,
\eeq
and the adjoint flux $\br:=-(A\nabla z+\bb z+\bg_2)\in\bSigma_d$ and $z$ solve
\beq \label{fos_d}
\divv\br+cz=g_1, \qquad \br+A\nabla z+\bb z=-\bg_2 \qquad \text{in }\O.
\eeq
By Green's formula, \eqref{fos_p} and \eqref{fos_d}, for every
$s\in H^1(\O)$,
\beq \label{green}
a(u,s)=\ell(s)-\langle\bsigma\cdot\bn,s\rangle_{\G_D},
\qquad
a(s,z)=m(s)-\langle\br\cdot\bn,s\rangle_{\G_D}.
\eeq

\subsection{Residuals, errors and least-squares estimators}

For $(\btau,w)$ and $(\brho,v)$ in $H(\divv;\O)\times H^1(\O)$, the
residuals of the constitutive and balance equations are
\begin{align} \label{res_def}
\Rcp(\btau,w)&:=\bff_2+\btau+A\nabla w, &
\Rbp(\btau,w)&:=\divv\btau+\bb\cdot\nabla w+cw-f_1, \\ \nonumber
\Rcd(\brho,v)&:=\bg_2+\brho+A\nabla v+\bb v, &
\Rbd(\brho,v)&:=\divv\brho+cv-g_1 .
\end{align}
They vanish at the exact solutions and are computable. For
$(\btau,w)\in H(\divv;\O)\times X_p$ and $(\brho,v)\in H(\divv;\O)\times X_d$
we write throughout
\beq \label{errors}
E_p:=\bsigma-\btau,\qquad e_p:=u-w,\qquad E_d:=\br-\brho,\qquad e_d:=z-v,
\eeq
so that $e_p,e_d\in H^1_D(\O)$, and $E_p,E_d\in H_N(\divv;\O)$ when
$\btau\in\bSigma_p$, $\brho\in\bSigma_d$. Then
\beq \label{res_err}
\begin{aligned}
\Rcp(\btau,w)&=-(E_p+A\nabla e_p), & \Rbp(\btau,w)&=-(\divv E_p+\bb\cdot\nabla e_p+ce_p),\\
\Rcd(\brho,v)&=-(E_d+A\nabla e_d+\bb e_d), & \Rbd(\brho,v)&=-(\divv E_d+ce_d).
\end{aligned}
\eeq
The least-squares functionals of \cite{WZ:26} and their estimators are
\begin{align} \label{LS}
\LS_p(\btau,w)&:=\|A^{-1/2}\Rcp(\btau,w)\|_0^2+\|\Rbp(\btau,w)\|_0^2=\eta_p^{\LS}(\btau,w)^2,
\\ \nonumber
\LS_d(\brho,v)&:=\|A^{-1/2}\Rcd(\brho,v)\|_0^2+\|\Rbd(\brho,v)\|_0^2=\eta_d^{\LS}(\brho,v)^2 .
\end{align}
Every bounded linear functional on $H^1_D(\O)$ can be written as
$(f_1,\cdot)-(\bff_2,\nabla\,\cdot)$. Since $\bff_2$ and $\bg_2$ enter the
constitutive equations, the residuals are in $L^2(\O)$ for such data
\cite{Zhang:23}. We write $\LS_p^0$ and $\LS_d^0$ for the functionals with homogeneous data,
$$
\LS_p^0(\bmu,s)=\|A^{-1/2}(\bmu+A\nabla s)\|_0^2
+\|\divv\bmu+\bb\cdot\nabla s+cs\|_0^2,
\qquad
\LS_d^0(\bmu,s)=\|A^{-1/2}(\bmu+A\nabla s+\bb s)\|_0^2+\|\divv\bmu+cs\|_0^2,
$$
so that, by \eqref{res_err}, $\LS_p(\btau,w)=\LS_p^0(E_p,e_p)$ and
$\LS_d(\brho,v)=\LS_d^0(E_d,e_d)$. We measure errors in the norm of
$H(\divv;\O)\times H^1(\O)$ weighted by $A$,
\beq \label{trinorm}
\tri(\bmu,s)\tri^2:=\|A^{-1/2}\bmu\|_0^2+\|\divv\bmu\|_0^2
+\|A^{1/2}\nabla s\|_0^2+\|s\|_0^2 .
\eeq

\begin{prop}[Norm equivalence]\label{prop:normeq}
There are constants $c_{\rm ls},C_{\rm ls}>0$ such that
$$
c_{\rm ls}\tri(\bmu,s)\tri\le\LS_p^0(\bmu,s)^{1/2}\le
C_{\rm ls}\tri(\bmu,s)\tri
\qquad\forall\,(\bmu,s)\in H_N(\divv;\O)\times H^1_D(\O),
$$
and likewise for $\LS_d^0$. Therefore, for arbitrary
$(\btau,w)\in\bSigma_p\times X_p$ and $(\brho,v)\in\bSigma_d\times X_d$,
\beq \label{normeq}
c_{\rm ls}\tri(E_p,e_p)\tri \le \eta_p^{\LS}(\btau,w)
\le C_{\rm ls}\tri(E_p,e_p)\tri, \qquad
c_{\rm ls}\tri(E_d,e_d)\tri \le \eta_d^{\LS}(\brho,v)
\le C_{\rm ls}\tri(E_d,e_d)\tri .
\eeq
\end{prop}

\noindent
In Proposition~\ref{prop:normeq}, the upper bound follows from the
triangle inequality and the boundedness of the coefficients, and the
lower bound is the coercivity of the least-squares functional of the
flux--potential system. This formulation appears already in
\cite{Jes:77}, and its coercivity in $H(\divv)\times H^1$ is classical
\cite{PCL:94,CLMM:94}. See also the monographs \cite{Jiang:98,BG:09} and
the survey \cite{BG:98}. For equations with convection and reaction terms
it is proved in \cite{CLMM:94}, in a weighted form for the adjoint system
in \cite[Thm.~2.2]{CFZ:15}, with the convection term partly replaced by
the flux in \cite[Thm.~4.1]{LZYZ:24}, and directly from the stability of
the second-order problem in \cite{Zhang:23}. The bounds \eqref{normeq}
follow from $\LS_p(\btau,w)=\LS_p^0(E_p,e_p)$ and
$\LS_d(\brho,v)=\LS_d^0(E_d,e_d)$. The least-squares functional already
serves as an a posteriori estimator in \cite{BMM:97}.

\begin{rem}[Weighting the balance term]\label{rem:weight}
For $\bb=\bzero$ and $c=1$ the equivalence holds with both constants equal
to one, because $(\bmu,\nabla s)=-(\divv\bmu,s)$, so that the cross terms
in $\LS_p^0$ cancel and $\LS_p^0=\LS_d^0=\tri\cdot\tri^2$. For
$\bb=\bzero$ and $c\ge c_0>0$, if we weight the balance terms of
\eqref{LS} by $c^{-1/2}$, the same computation gives
$$
\|A^{-1/2}(\bmu+A\nabla s)\|_0^2+\|c^{-1/2}(\divv\bmu+cs)\|_0^2
=\|A^{-1/2}\bmu\|_0^2+\|c^{-1/2}\divv\bmu\|_0^2
+\|A^{1/2}\nabla s\|_0^2+\|c^{1/2}s\|_0^2
$$
for every such $c$, so the weighted estimators equal the errors in the
weighted norm. This is the weighting of the robust estimator of
\cite{AV:14} for reaction--diffusion problems. The same weighting arises
when we write these problems as convex minimization problems
\cite{Zhang:20}, since the weighted functional is twice the primal--dual
gap.
\end{rem}

\subsection{The evaluators \texorpdfstring{$\gamma_0$}{gamma0} and \texorpdfstring{$\gamma_1$}{gamma1}}

Following Giles and S\"uli \cite{GS:02}, let
\beq \label{gamma0}
\gamma_0(w,v):=m(w)+\ell(v)-a(w,v), \qquad w,v\in H^1(\O).
\eeq
Then
\beq \label{gap_gamma0}
\gamma-\gamma_0(w,v)=a(e_p,e_d)
\eeq
for all $w\in X_p$, $v\in X_d$, with no orthogonality assumption
\cite{WZ:26}.

If fluxes are available, we can express the gradients in \eqref{gamma0}
through them, $\nabla u=-A^{-1}(\bsigma+\bff_2)$ and
$A\nabla z+\bb z=-(\br+\bg_2)$. Substituting these expressions in $m$,
$\ell$ and $a$, we obtain $\tilde m$, $\tilde\ell$ and $\tilde a$ below. The
second evaluator of \cite{WZ:26} adds to $\tilde m+\tilde\ell-\tilde a$
two corrections, the constitutive residuals tested with the fluxes, which
vanish at the exact solutions:
\begin{align} \label{gamma1}
\gamma_1(w,\btau,v,\brho)
&:= \tilde m(\btau,w)+\tilde\ell(\brho,v)
-\tilde a((\btau,w),(\brho,v))
+ \big(\Rcp(\btau,w),\,A^{-1}\brho\big) + \big(A^{-1}\btau,\,\Rcd(\brho,v)\big),
\\ \nonumber
\tilde m(\btau,w)&:=(g_1,w)+(A^{-1}(\btau+\bff_2),\bg_2)
-\langle\psi_d,w\rangle_{\G_N},
\\ \nonumber
\tilde\ell(\brho,v)&:=(f_1,v)+(A^{-1}(\brho+\bb v+\bg_2),\bff_2)
-\langle\psi_p,v\rangle_{\G_N},
\\ \nonumber
\tilde a((\btau,w),(\brho,v))&:=(A^{-1}(\btau+\bff_2),\brho+\bg_2)+(cw,v),
\end{align}
for $(\btau,w),(\brho,v)\in H(\divv;\O)\times H^1(\O)$. We need the
boundary conditions $w\in X_p$ and $v\in X_d$ only for consistency, as
in \eqref{gap_gamma0}. At the exact solutions $\Rcp(\bsigma,u)=\Rcd(\br,z)=0$,
and $\gamma_1(u,\bsigma,z,\br)=\gamma_0(u,z)=\gamma$.

\section{The QoI gaps between the evaluators}\label{sec:gap}

\subsection{The QoI gap identity}\label{sec:identity}

\begin{thm}[QoI gap identity]\label{thm:gap}
For all $(\btau,w)$ and $(\brho,v)$ in $H(\divv;\O)\times H^1(\O)$,
\beq \label{gap}
\gamma_1(w,\btau,v,\brho)-\gamma_0(w,v)
= \big(A^{-1}\Rcp(\btau,w),\;\Rcd(\brho,v)\big).
\eeq
\end{thm}

\begin{proof}
By the definitions of $\tilde m$ and $m$, and of $\tilde\ell$ and $\ell$,
$$
\tilde m(\btau,w)=m(w)+\big(\Rcp(\btau,w),A^{-1}\bg_2\big),\qquad
\tilde\ell(\brho,v)=\ell(v)+\big(\Rcd(\brho,v),A^{-1}\bff_2\big).
$$
Since $A$ is symmetric,
$a(w,v)-(cw,v)=\big(A\nabla w,\,A^{-1}(A\nabla v+\bb v)\big)$.
We subtract \eqref{gamma0} from \eqref{gamma1}. The terms $m(w)$ and
$\ell(v)$ cancel, and
\begin{align*}
&\gamma_1(w,\btau,v,\brho)-\gamma_0(w,v)\\
&\quad=\underbrace{\big(\Rcp(\btau,w),A^{-1}\bg_2\big)+\big(\Rcp(\btau,w),A^{-1}\brho\big)}
_{=\,(\Rcp(\btau,w),\,A^{-1}(\bg_2+\brho))}
+\underbrace{\big(\Rcd(\brho,v),A^{-1}\bff_2\big)+\big(A^{-1}\btau,\Rcd(\brho,v)\big)}
_{=\,(A^{-1}(\bff_2+\btau),\,\Rcd(\brho,v))} \\
&\qquad -\big(A^{-1}(\btau+\bff_2),\,\brho+\bg_2\big)
+\big(A\nabla w,\,A^{-1}(A\nabla v+\bb v)\big).
\end{align*}
We split $\Rcd(\brho,v)=(\bg_2+\brho)+(A\nabla v+\bb v)$ in the second bracket.
Its first part cancels the third term, and we are left with
\begin{align*}
\gamma_1(w,\btau,v,\brho)-\gamma_0(w,v)
&=\big(\Rcp(\btau,w),\,A^{-1}(\bg_2+\brho)\big)
+\big(A^{-1}(\bff_2+\btau+A\nabla w),\,A\nabla v+\bb v\big)\\
&=\big(A^{-1}\Rcp(\btau,w),\,\Rcd(\brho,v)\big). \qedhere
\end{align*}
\end{proof}

\noindent
The proof uses only $L^2$ inner products, so the identity also holds for
fluxes that are only in $L^2(\O)^d$, for instance $-(A\nabla\tilde w+\bff_2)$
for another potential $\tilde w$, or the broken flux of a
Crouzeix--Raviart approximation. We take all fluxes in $H(\divv;\O)$
because the balance residuals need it, and so do the least-squares
bounds and $\gamma_2(\btau,\brho)$, which use these residuals. The QoI
gap itself uses only the
$L^2$ errors of the fluxes.

\begin{cor}\label{cor:cs}
For $(\btau,w)\in\bSigma_p\times X_p$ and $(\brho,v)\in\bSigma_d\times X_d$,
\beq \label{csbound}
\begin{aligned}
|\gamma_1(w,\btau,v,\brho)-\gamma_0(w,v)|
&\le \|A^{-1/2}\Rcp(\btau,w)\|_0\,\|A^{-1/2}\Rcd(\brho,v)\|_0\\
&\le \eta_p^{\LS}(\btau,w)\,\eta_d^{\LS}(\brho,v).
\end{aligned}
\eeq
\end{cor}
\begin{proof}
By Theorem~\ref{thm:gap}, $\gamma_1(w,\btau,v,\brho)-\gamma_0(w,v)$ is the
inner product of $A^{-1/2}\Rcp(\btau,w)$ and $A^{-1/2}\Rcd(\brho,v)$.
The first inequality is the Cauchy--Schwarz inequality, and the second
holds because each factor is a part of the corresponding estimator
\eqref{LS}.
\end{proof}

We call the computable quantity
\beq \label{delta}
\begin{aligned}
\delta(w,\btau,v,\brho) &:= \gamma_1(w,\btau,v,\brho)-\gamma_0(w,v)\\
&\phantom{:}= \big(A^{-1}(\bff_2+\btau+A\nabla w),\;\bg_2+\brho+A\nabla v+\bb v\big)
\end{aligned}
\eeq
the {\em signed QoI gap}, or simply the {\em QoI gap}. By \eqref{res_err},
\beq \label{repr}
\delta(w,\btau,v,\brho)=\big(A^{-1}(E_p+A\nabla e_p),\,E_d+A\nabla e_d+\bb e_d\big)
\eeq
for $(\btau,w)\in H(\divv;\O)\times X_p$ and $(\brho,v)\in H(\divv;\O)\times X_d$.
The QoI gap is thus the $A^{-1}$ inner product of the errors of the two
constitutive equations.
The error of the potential evaluator,
$$
\gamma-\gamma_0(w,v)=a(e_p,e_d)=(\nabla e_p,A\nabla e_d+\bb e_d)+(ce_p,e_d),
$$
is the $A^{-1}$ inner product of the potential parts of the two constitutive errors plus
the reaction term. Since
$$
\gamma-\gamma_0(w,v)=\delta(w,\btau,v,\brho)+\big(\gamma-\gamma_1(w,\btau,v,\brho)\big),
$$
the computable QoI gap equals the error of $\gamma_0(w,v)$ up to the error of
$\gamma_1(w,\btau,v,\brho)$.

\subsection{Bounds and two-sided estimates}\label{sec:estimate}

\begin{cor}[Bounds]\label{cor:bounds}
For $(\btau,w)\in\bSigma_p\times X_p$ and $(\brho,v)\in\bSigma_d\times X_d$,
\beq \label{bound_gamma1}
\begin{aligned}
|\gamma-\gamma_0(w,v)|&\le C\,\eta_p^{\LS}(\btau,w)\,\eta_d^{\LS}(\brho,v),\\
|\gamma-\gamma_1(w,\btau,v,\brho)|&\le C\,\eta_p^{\LS}(\btau,w)\,\eta_d^{\LS}(\brho,v).
\end{aligned}
\eeq
\end{cor}
\begin{proof}
We bound \eqref{gap_gamma0} by the boundedness of $a$ and use
\eqref{normeq}. For $\gamma_1(w,\btau,v,\brho)$ we use
$\gamma-\gamma_1(w,\btau,v,\brho)=\gamma-\gamma_0(w,v)-\delta(w,\btau,v,\brho)$
and \eqref{csbound}.
\end{proof}

Whether a QoI gap measures the errors of the two evaluators depends on the
signs of the errors. The argument uses only that a QoI gap is the difference
of two approximations of $\gamma$, so we state it for two arbitrary
approximations.

\begin{prop}[Two-sided estimate]\label{prop:twosided}
Let $\gamma_A$ and $\gamma_B$ be two approximations of $\gamma$ with
errors $e_A:=\gamma-\gamma_A$ and $e_B:=\gamma-\gamma_B$, labeled so that
$|e_A|\le|e_B|$. Then $|\gamma_A-\gamma_B|\le|e_A|+|e_B|\le2|e_B|$, and
the following hold.
\begin{enumerate}
\item[(a)] If $e_Ae_B\le0$, that is, if $\gamma$ lies between $\gamma_A$
and $\gamma_B$, then $|\gamma_A-\gamma_B|=|e_A|+|e_B|$, and consequently
$|e_B|\le|\gamma_A-\gamma_B|\le2|e_B|$.
\item[(b)] If $e_Ae_B>0$, then
$|\gamma_A-\gamma_B|=|e_B|-|e_A|=\big(1-|e_A|/|e_B|\big)|e_B|$.
\end{enumerate}
In particular, $\gamma$ lies between $\gamma_A$ and $\gamma_B$ if and only
if $|\gamma_A-\gamma_B|=|e_A|+|e_B|$.
\end{prop}
\begin{proof}
If $e_Ae_B\le0$ then $|e_B-e_A|=|e_A|+|e_B|$, and if $e_Ae_B>0$ then
$|e_B-e_A|=|e_B|-|e_A|$, which equals $|e_A|+|e_B|$ only if $e_A=0$,
excluded in case (b). The remaining statements follow from
$|e_A|\le|e_B|$.
\end{proof}

\noindent
We apply the proposition to $\gamma_0(w,v)$ and
$\gamma_1(w,\btau,v,\brho)$, with $|\gamma_A-\gamma_B|=|\delta(w,\btau,v,\brho)|$.
The QoI gap then estimates the larger of the two errors. It does so
within a factor two when $\gamma$ lies between the two evaluators, and
within the factor $(1-|e_A|/|e_B|)^{-1}$ otherwise. This factor is large
only when the two evaluators make nearly the same error, in sign and in
size. In either case we do not need to know which evaluator is the more
accurate. When one
error is much smaller than the other, the QoI gap is a signed estimate of
the error of the less accurate evaluator, since
$\gamma_A-\gamma_B=e_B-e_A\approx e_B$.

The value that we report for the quantity of interest need not be one of
the evaluators. In case (a), $\gamma$ lies in the interval between them,
so the length of the interval bounds the error of any value in it, and
half of the length bounds the error of the midpoint. These errors can be
much smaller than the QoI gap.

\begin{rem}[When a QoI gap vanishes]\label{rem:zero}
If a constitutive residual vanishes identically although the pair is not
exact, that is, $\btau=-(A\nabla w+\bff_2)$ or
$\brho=-(A\nabla v+\bb v+\bg_2)$, then the flux is the flux of the
potential, $\gamma_1(w,\btau,v,\brho)=\gamma_0(w,v)$, and the QoI gap is zero
whatever the errors are. Theorem~\ref{thm:gap} allows such a flux
whenever $A\nabla w+\bff_2\in H(\divv;\O)$. So we must compute the flux
independently, by a least-squares method or a recovery, and not by
differentiating the potential.

A QoI gap also vanishes when both errors are nonzero but equal, in case
(b) of Proposition~\ref{prop:twosided}. Along a refinement sequence on which the two
errors keep the same sign and the more accurate of the two evaluators
changes, the difference of the two absolute errors changes sign, and
with it the QoI gap.
\end{rem}

\subsection{The reaction-restored evaluator \texorpdfstring{$\gamma_2$}{gamma2}}\label{sec:gamma2}

The signed QoI gap does not contain the reaction term $(ce_p,e_d)$ of
$\gamma-\gamma_0(w,v)$. The coefficient $c$ occurs only in the balance
equations of \eqref{fos_p} and \eqref{fos_d},
while the signed QoI gap pairs the constitutive residuals, so $c$ enters it
only through the computed fields. Since
$|(c\,e_p,e_d)|\le\|c\|_{0,\infty}\|e_p\|_0\|e_d\|_0$, this term is of
higher order than the product of the $H^1$ errors when $c$ is not large
and the $L^2$ errors of the potentials are of higher order than their
$H^1$ errors. For large $c$ the factor $\|c\|_{0,\infty}$ can make the
term dominant, and the reaction term causes layers that the mesh may not
resolve.

To restore the reaction term, we solve the balance equations for the
potentials, as in the flux-only least-squares method for linear
transport equations \cite{LZ:19}. We express the convection term of the primal balance
equation through the flux, $\bb\cdot\nabla u=-\bb\cdot A^{-1}(\bsigma+\bff_2)$,
as in mixed methods, so that the balance equation becomes
$\divv\bsigma-\bb\cdot A^{-1}(\bsigma+\bff_2)+cu=f_1$, one of the
first-order systems of \cite[\S4]{LZYZ:24}. Its residual is
$\Rbp(\btau,w)-\bb\cdot A^{-1}\Rcp(\btau,w)$. The adjoint balance equation of \eqref{fos_d}
contains no convection term. For $c\ge c_0>0$ we solve these equations
for the potentials: for $\btau\in\bSigma_p$ and $\brho\in\bSigma_d$
let
\beq \label{whatvhat}
\hat w(\btau):=c^{-1}\big(f_1-\divv\btau+\bb\cdot A^{-1}(\btau+\bff_2)\big),
\qquad
\hat v(\brho):=c^{-1}(g_1-\divv\brho),
\eeq
so that $u=\hat w(\bsigma)$ and $z=\hat v(\br)$.

By \eqref{green} and \eqref{QoI},
$\gamma=a(u,z)+\langle\bsigma\cdot\bn,\phi_d\rangle_{\G_D}
+\langle\br\cdot\bn,\phi_p\rangle_{\G_D}$, with
$a(u,z)=(\nabla u,A\nabla z+\bb z)+(cu,z)$. The constitutive equations,
$\nabla u=-A^{-1}(\bsigma+\bff_2)$ and $A\nabla z+\bb z=-(\br+\bg_2)$, give
$$
\gamma=\big(A^{-1}(\bsigma+\bff_2),\,\br+\bg_2\big)+(cu,z)
+\langle\bsigma\cdot\bn,\phi_d\rangle_{\G_D}+\langle\br\cdot\bn,\phi_p\rangle_{\G_D}.
$$
For $c=0$ there is no reaction term, and this formula expresses $\gamma$
through the fluxes alone. For $c\ge c_0>0$, replacing $u$ and $z$
by $\hat w(\bsigma)$ and $\hat v(\br)$ gives $\gamma=\gamma_2(\bsigma,\br)$,
where the reaction-restored evaluator
\beq \label{gamma2}
\gamma_2(\btau,\brho):=\big(A^{-1}(\btau+\bff_2),\,\brho+\bg_2\big)
+\big(c\,\hat w(\btau),\hat v(\brho)\big)
+\langle\btau\cdot\bn,\phi_d\rangle_{\G_D}+\langle\brho\cdot\bn,\phi_p\rangle_{\G_D}
\eeq
is a functional of the two fluxes alone, just as \eqref{gamma0} uses
only the potentials. The evaluator $\gamma_0$ pairs each Neumann datum
with the trace of the other potential, and the potentials satisfy the
Dirichlet conditions. The evaluator $\gamma_2$ pairs each Dirichlet datum
with the normal trace of the other flux, and the fluxes satisfy the
Neumann conditions. For potentials $w\in X_p$ and $v\in X_d$,
$\gamma_1(w,\btau,v,\brho)$ has the same form. It differs from
$\gamma_2(\btau,\brho)$ only in the reaction term, which it evaluates with
both the potentials $w$, $v$ and the potentials $\hat w(\btau)$,
$\hat v(\brho)$ obtained from the fluxes:
\beq \label{gamma1flux}
\begin{aligned}
\gamma_1(w,\btau,v,\brho)&=\big(A^{-1}(\btau+\bff_2),\,\brho+\bg_2\big)
+(c\hat w(\btau),v)+(cw,\hat v(\brho))-(cw,v)\\
&\quad+\langle\btau\cdot\bn,\phi_d\rangle_{\G_D}+\langle\brho\cdot\bn,\phi_p\rangle_{\G_D},\\
\gamma_2(\btau,\brho)-\gamma_1(w,\btau,v,\brho)&=\big(c(\hat w(\btau)-w),\,\hat v(\brho)-v\big).
\end{aligned}
\eeq
The second identity pairs the differences between the potentials
obtained from the fluxes and the potentials $w$ and $v$.

\begin{prop}[Reaction-restored evaluator]\label{prop:gamma2}
Let $c\ge c_0>0$, $\btau\in\bSigma_p$ and $\brho\in\bSigma_d$. Then
\beq \label{err_gamma2}
\gamma-\gamma_2(\btau,\brho)=-\big(A^{-1}E_p,E_d\big)
-\big(c^{-1}(\divv E_p-\bb\cdot A^{-1}E_p),\,\divv E_d\big),
\eeq
and for all $w\in X_p$ and $v\in X_d$, \eqref{gamma1flux} holds,
\begin{align}
\gamma_2(\btau,\brho)&=\gamma_1(w,\btau,v,\brho)
+\big(c^{-1}(\Rbp(\btau,w)-\bb\cdot A^{-1}\Rcp(\btau,w)),\,\Rbd(\brho,v)\big)
\notag\\
&=\gamma_0(w,v)+\big(A^{-1}\Rcp(\btau,w),\Rcd(\brho,v)\big)
+\big(c^{-1}\Rbp(\btau,w),\Rbd(\brho,v)\big)\notag\\
&\quad-\big(c^{-1}\bb\cdot A^{-1}\Rcp(\btau,w),\,\Rbd(\brho,v)\big),
\label{gamma2res}\\
\gamma_2(\btau,\brho)-\gamma_0(w,v)&=a(e_p,e_d)+\big(A^{-1}E_p,E_d\big)
+\big(c^{-1}(\divv E_p-\bb\cdot A^{-1}E_p),\,\divv E_d\big).
\label{delta2}
\end{align}
\end{prop}
\begin{proof}
We insert $\btau=\bsigma-E_p$ and $\brho=\br-E_d$ in \eqref{gamma2}, so
that $c\hat w(\btau)=cu+\divv E_p-\bb\cdot A^{-1}E_p$ and
$c\hat v(\brho)=cz+\divv E_d$, and use $\gamma=\gamma_2(\bsigma,\br)$ and
the constitutive equations. The terms linear in $E_p$ are
$(E_p,\nabla z)+(\divv E_p,z)-\langle E_p\cdot\bn,\phi_d\rangle_{\G_D}$, and
those linear in $E_d$ are
$(\nabla u,E_d)+(u,\divv E_d)-\langle E_d\cdot\bn,\phi_p\rangle_{\G_D}$.
Both vanish by Green's formula, since $E_p\cdot\bn=E_d\cdot\bn=0$ on
$\G_N$, $z=\phi_d$ and $u=\phi_p$ on $\G_D$. The remaining terms are
quadratic and give \eqref{err_gamma2}.

By \eqref{gap} and \eqref{gamma0},
$\gamma_1(w,\btau,v,\brho)=m(w)+\ell(v)-a(w,v)+(A^{-1}\Rcp(\btau,w),\Rcd(\brho,v))$. We expand
$(A^{-1}\Rcp(\btau,w),\Rcd(\brho,v))$ by \eqref{res_def}. The terms $(\bg_2,\nabla w)$,
$(\bff_2,\nabla v)$ and $(A\nabla w,\nabla v)+(\bb\cdot\nabla w,v)$ cancel
against those of $m(w)$, $\ell(v)$ and $a(w,v)$, and Green's formula
applied to $(\btau,\nabla v)$ and $(\nabla w,\brho)$, with the boundary
conditions of $w$, $v$, $\btau$ and $\brho$, gives \eqref{gamma1flux}.
Since $c(\hat w(\btau)-w)=-(\Rbp(\btau,w)-\bb\cdot A^{-1}\Rcp(\btau,w))$ and
$\hat v(\brho)-v=-c^{-1}\Rbd(\brho,v)$, this gives the first equality in
\eqref{gamma2res}, and \eqref{gap} gives the second. Finally, we obtain \eqref{delta2} by subtracting \eqref{err_gamma2}
from \eqref{gap_gamma0}.
\end{proof}

\begin{rem}[The case $c=0$]\label{rem:c0}
For $c=0$ the potentials \eqref{whatvhat} are not defined, but
$c\hat w(\btau)=f_1-\divv\btau+\bb\cdot A^{-1}(\btau+\bff_2)$ and
$c\hat v(\brho)=g_1-\divv\brho$ are, and the proof of \eqref{gamma1flux} applies
unchanged:
$$
\begin{aligned}
\gamma_1(w,\btau,v,\brho)&=\big(A^{-1}(\btau+\bff_2),\,\brho+\bg_2\big)
+\big(f_1-\divv\btau+\bb\cdot A^{-1}(\btau+\bff_2),\,v\big)
+\big(w,\,g_1-\divv\brho\big)\\
&\quad+\langle\btau\cdot\bn,\phi_d\rangle_{\G_D}+\langle\brho\cdot\bn,\phi_p\rangle_{\G_D},
\end{aligned}
$$
where the second and third terms test the balance residuals of the
fluxes with the potentials. For fluxes that satisfy the balance
equations, these two terms vanish. Then
$\gamma_1(w,\btau,v,\brho)$ does not depend on $w$ and $v$, as
$\gamma_0(w,v)$ does not depend on the fluxes, and, as in the proof of
Proposition~\ref{prop:gamma2},
$\gamma-\gamma_1(w,\btau,v,\brho)=-(A^{-1}E_p,E_d)$. For $\bb=\bzero$ these are the
equilibrated fluxes. The fluxes of a mixed method
are such fluxes, up to the oscillation of $f_1$ and $g_1$. Then $\gamma_1(w,\btau,v,\brho)$
evaluates $\gamma$ from the mixed fluxes alone, and the QoI gap compares it
with $\gamma_0(w,v)$, which uses only the potentials. If the balance residuals
of the fluxes are of order $c$, the reaction term $(c\hat w(\btau),\hat v(\brho))$ of
\eqref{gamma2} tends to zero as $c\to0$, and $\gamma_2(\btau,\brho)$ tends to the same
expression.
\end{rem}

\noindent
The error \eqref{err_gamma2} of $\gamma_2(\btau,\brho)$ pairs the primal
and adjoint flux errors, as the error $a(e_p,e_d)$ of $\gamma_0(w,v)$
pairs the potential errors. We call
$$
\delta_2(w,\btau,v,\brho):=\gamma_2(\btau,\brho)-\gamma_0(w,v)
$$
the {\em reaction-restored QoI gap}. By \eqref{delta2} it equals the whole error of the
potential evaluator, reaction term included, up to the error of
$\gamma_2(\btau,\brho)$. The three evaluators use more and more of the
fluxes. By \eqref{gap} and \eqref{gamma2res},
$\gamma_1(w,\btau,v,\brho)-\gamma_0(w,v)$ pairs the constitutive
residuals, $\gamma_2(\btau,\brho)-\gamma_1(w,\btau,v,\brho)$ pairs the
balance residuals, with a convective coupling term, and
$\delta_2(w,\btau,v,\brho)$ pairs both. We measure both QoI gaps from the
potential evaluator $\gamma_0(w,v)$. Proposition~\ref{prop:twosided}
applies to $\delta_2(w,\btau,v,\brho)$ in the same way, with
$\gamma_0(w,v)$ and $\gamma_2(\btau,\brho)$ as the two approximations,
and also to $\gamma_1(w,\btau,v,\brho)$ and $\gamma_2(\btau,\brho)$. We
bound the reaction-restored QoI gap by the least-squares estimators with the balance
terms weighted by $c^{-1/2}$, as in Remark~\ref{rem:weight},
\begin{align} \label{LSc}
\eta_{p,c}^{\LS}(\btau,w)&:=\big(\|A^{-1/2}\Rcp(\btau,w)\|_0^2+\|c^{-1/2}\Rbp(\btau,w)\|_0^2\big)^{1/2},
\\ \nonumber
\eta_{d,c}^{\LS}(\brho,v)&:=\big(\|A^{-1/2}\Rcd(\brho,v)\|_0^2+\|c^{-1/2}\Rbd(\brho,v)\|_0^2\big)^{1/2}.
\end{align}

\begin{prop}[Reaction-restored QoI gap]\label{prop:delta2}
Let $c\ge c_0>0$, let $(\btau,w)\in\bSigma_p\times X_p$ and
$(\brho,v)\in\bSigma_d\times X_d$, and let
$\vartheta:=\|c^{-1/2}A^{-1/2}\bb\|_{0,\infty}$. Then
\beq \label{delta2res}
\begin{aligned}
\delta_2(w,\btau,v,\brho)&=\big(A^{-1}\Rcp(\btau,w),\Rcd(\brho,v)\big)
+\big(c^{-1}\Rbp(\btau,w),\Rbd(\brho,v)\big)\\
&\quad-\big(c^{-1}\bb\cdot A^{-1}\Rcp(\btau,w),\,\Rbd(\brho,v)\big)
\end{aligned}
\eeq
and
\beq \label{delta2bound}
|\delta_2(w,\btau,v,\brho)|\le(1+\vartheta)\,\eta_{p,c}^{\LS}(\btau,w)\,\eta_{d,c}^{\LS}(\brho,v).
\eeq
For $\bb=\bzero$, $\delta_2(w,\btau,v,\brho)$ is the inner product in
$L^2(\O)^d\times L^2(\O)$, with the weights $A^{-1}$ and $c^{-1}$, of the
residual pairs of the two problems. The norm induced by this inner
product is the weighted estimator \eqref{LSc}, so
$|\delta_2(w,\btau,v,\brho)|\le\eta_{p,c}^{\LS}(\btau,w)\,\eta_{d,c}^{\LS}(\brho,v)$.
\end{prop}
\begin{proof}
The representation \eqref{delta2res} is \eqref{gamma2res}, by the
definition of $\delta_2(w,\btau,v,\brho)$. By the Cauchy--Schwarz inequality in
$L^2(\O)$ and then in $\mathbb R^2$,
$$
\begin{aligned}
&\big|(A^{-1}\Rcp(\btau,w),\Rcd(\brho,v))+(c^{-1}\Rbp(\btau,w),\Rbd(\brho,v))\big|\\
&\quad\le\|A^{-1/2}\Rcp(\btau,w)\|_0\|A^{-1/2}\Rcd(\brho,v)\|_0
+\|c^{-1/2}\Rbp(\btau,w)\|_0\|c^{-1/2}\Rbd(\brho,v)\|_0\\
&\quad\le\eta_{p,c}^{\LS}(\btau,w)\,\eta_{d,c}^{\LS}(\brho,v).
\end{aligned}
$$
Pointwise,
$$
|c^{-1/2}\bb\cdot A^{-1}\Rcp(\btau,w)|
=c^{-1/2}|(A^{-1/2}\bb)\cdot(A^{-1/2}\Rcp(\btau,w))|\le\vartheta|A^{-1/2}\Rcp(\btau,w)|,
$$
so that
$$
\begin{aligned}
|(c^{-1}\bb\cdot A^{-1}\Rcp(\btau,w),\Rbd(\brho,v))|
&\le\vartheta\|A^{-1/2}\Rcp(\btau,w)\|_0\|c^{-1/2}\Rbd(\brho,v)\|_0\\
&\le\vartheta\,\eta_{p,c}^{\LS}(\btau,w)\,\eta_{d,c}^{\LS}(\brho,v).
\end{aligned}
$$
\end{proof}

\begin{cor}\label{cor:bound2}
Let $c\ge c_0>0$, $(\btau,w)\in\bSigma_p\times X_p$ and
$(\brho,v)\in\bSigma_d\times X_d$. Then
\beq \label{bound_gamma2}
|\gamma-\gamma_2(\btau,\brho)|\le C\,\eta_{p,c}^{\LS}(\btau,w)\,\eta_{d,c}^{\LS}(\brho,v).
\eeq
\end{cor}
\begin{proof}
We bound each term of \eqref{err_gamma2} by the Cauchy--Schwarz
inequality, using $|c^{-1/2}\bb\cdot A^{-1}E_p|\le\vartheta|A^{-1/2}E_p|$,
and use \eqref{normeq}. Since $c_0\le c\le\|c\|_{0,\infty}$,
$$
\min\{1,\|c\|_{0,\infty}^{-1}\}\,\eta_p^{\LS}(\btau,w)^2\le\eta_{p,c}^{\LS}(\btau,w)^2
\le\max\{1,c_0^{-1}\}\,\eta_p^{\LS}(\btau,w)^2,
$$
and likewise for the adjoint estimators.
\end{proof}

\begin{rem}[The errors of the three evaluators]\label{rem:gamma2}
The error $a(e_p,e_d)$ of $\gamma_0(w,v)$ contains the reaction term
$(ce_p,e_d)$, and so does the error of
$\gamma_1(w,\btau,v,\brho)=\gamma_0(w,v)+\delta(w,\btau,v,\brho)$, since the
signed QoI gap does not contain it. The error \eqref{err_gamma2} of
$\gamma_2(\btau,\brho)$ does not contain it. Instead, it contains the
inner product of the two divergence errors weighted by $c^{-1}$ and the
cross term $(c^{-1}\bb\cdot A^{-1}E_p,\divv E_d)$, which is bounded by
$\vartheta\|A^{-1/2}E_p\|_0\|c^{-1/2}\divv E_d\|_0$ with $\vartheta$ of
Proposition~\ref{prop:delta2}. For $A=\epsilon I$ with constant
$\epsilon$, $\bb$ and $c$, we have $\vartheta=|\bb|/\sqrt{\epsilon c}$,
which is large when $c$ is small compared with $|\bb|^2/\epsilon$.
Where reaction dominates, the reaction term can dominate the errors of
$\gamma_0(w,v)$ and $\gamma_1(w,\btau,v,\brho)$, while the error of
$\gamma_2(\btau,\brho)$ does not contain it. However, the terms
weighted by $c^{-1}$ in \eqref{err_gamma2} can make $\gamma_2(\btau,\brho)$
less accurate than $\gamma_1(w,\btau,v,\brho)$, especially through the
cross term when $\vartheta$ is large. We do not use $\gamma_2(\btau,\brho)$
there. Unlike the signed QoI gap (Remark~\ref{rem:zero}), the
reaction-restored QoI gap does not vanish when a flux is the flux of its
potential. If $\btau=-(A\nabla w+\bff_2)$, then
$\delta_2(w,\btau,v,\brho)=(c^{-1}\Rbp(\btau,w),\Rbd(\brho,v))$, which
still pairs the primal and adjoint balance residuals.
\end{rem}

\section{Relation to other approaches}\label{sec:relation}

A posteriori error estimation rests on two principles. The first
compares two approximations of the same quantity, as when we compare a
quadrature rule on a mesh with the same rule on its refinement, and gives
an estimate. The second evaluates the residual, as we do for a linear
system $Ax=b$ through $A(x-x_h)=b-Ax_h$, and gives a reliable bound,
possibly with an overestimating constant \cite{LZYZ:24}. The
dual-weighted residual estimate is an instance of the first principle and
the product bound of the second.
Recovery-based estimators belong to the first, since they compare the
flux of the numerical potential with a recovered flux, and residual-based
estimators to the second. The least-squares functional contains both. Its
constitutive residual $\Rcp(\btau,w)$ is the difference of the computed flux
$\btau$ and the flux $-(A\nabla w+\bff_2)$ of the potential, as in the
first principle. Its balance residual $\Rbp(\btau,w)$ is a residual in the
sense of the second. Together they are reliable for arbitrary arguments
(Proposition~\ref{prop:normeq}). In this paper we
apply both principles to the same residuals. The QoI gap compares the two
approximations $\gamma_0(w,v)$ and $\gamma_1(w,\btau,v,\brho)$ of $\gamma$
without assuming that either of them is the more accurate, and it
estimates the larger of their errors. The norms of the residuals, the
least-squares estimators, give the reliable upper bound, and we use them
to refine the mesh.

\subsection{The dual-weighted residual estimate and the QoI gap}\label{sec:dwr}

We build the dual-weighted residual estimate and the QoI gap in the same way:
we compute $\gamma$ twice and take the difference of the two values as
the estimate. The dual-weighted residual method evaluates $\gamma_0$
twice, with the same potential $w$ and two adjoint approximations, $v$
and a more accurate $\tilde v\in X_d$, computed in an enriched space or
by patchwise higher-order interpolation of $v$ \cite[\S5.1]{BR:01}:
\beq \label{etadwr}
\eta_{\rm DWR}:=\gamma_0(w,\tilde v)-\gamma_0(w,v)
=\ell(\tilde v-v)-a(w,\tilde v-v)=a(u-w,\tilde v-v),
\eeq
that is, $\gamma-\gamma_0(w,v)=a(u-w,z-v)$ with $z-v$ replaced by
$\tilde v-v$. If $w$
is a Galerkin approximation and $v$ lies in its test space, as for a
Galerkin approximation $v$ on the same mesh with $\phi_d=0$, Galerkin
orthogonality gives $\gamma_0(w,v)=m(w)$, so that \eqref{etadwr} estimates the error of the
direct output.

Since the estimate is the difference of two values of $\gamma$, other
choices of the two values give variants of it. We do not need Galerkin
orthogonality. The identity \eqref{gap_gamma0} holds for all $w\in X_p$
and $v\in X_d$, so \eqref{etadwr} holds for arbitrary approximations and then
concerns the error of $\gamma_0(w,v)$. This value is the output $m(w)$
corrected by the primal residual tested with $v$ \cite{GS:02}, and the
orthogonality only makes this correction vanish. We can also change the
potentials. With a more accurate primal approximation $\tilde w\in X_p$,
the difference $\gamma_0(\tilde w,v)-\gamma_0(w,v)=a(\tilde w-w,z-v)$
replaces $u-w$ by $\tilde w-w$, and we can average the two one-sided
differences. With both,
$\gamma_0(\tilde w,\tilde v)-\gamma_0(w,v)
=\gamma-\gamma_0(w,v)-a(u-\tilde w,z-\tilde v)$. Finally, we can keep the
potentials and use fluxes instead of more accurate potentials. If we
replace the gradients $\nabla w$ and $A\nabla v+\bb v$ in $\gamma_0(w,v)$
by $-A^{-1}(\btau+\bff_2)$ and $-(\brho+\bg_2)$ and add the two residual
corrections of \eqref{gamma1}, we obtain $\gamma_1(w,\btau,v,\brho)$, and
the difference of the two values is the signed QoI gap
$\delta(w,\btau,v,\brho)$.

The dual-weighted residual estimate and the QoI gap are both differences of two values of $\gamma$, so
Proposition~\ref{prop:twosided} applies to both. The two
estimates differ in what they assume. The dual-weighted residual estimate
$\eta_{\rm DWR}=\gamma-\gamma_0(w,v)-a(u-w,z-\tilde v)$ estimates $\gamma-\gamma_0(w,v)$ when $\gamma_0(w,\tilde v)$ is much more accurate than
$\gamma_0(w,v)$, while the QoI gap does not assume that either of its two
values is the more accurate.

The QoI gap is a generalized dual-weighted residual estimate with
flux-based weights. In $\gamma-\gamma_0(w,v)=(\nabla e_p,A\nabla e_d+\bb
e_d)+(ce_p,e_d)$, the dual-weighted residual estimate replaces $z-v$ by
$\tilde v-v$. The QoI gap replaces $\nabla e_p$ by
$-A^{-1}\Rcp(\btau,w)=\nabla e_p+A^{-1}E_p$ and $A\nabla e_d+\bb e_d$ by
$-\Rcd(\brho,v)=A\nabla e_d+\bb e_d+E_d$, which are exact when the fluxes
are exact. It also omits the reaction term, and the result is
$\delta(w,\btau,v,\brho)$. The weights come from the fluxes, on both the
primal and the adjoint side, and the fluxes need not be more accurate than
the gradients of the potentials.

\subsection{The QoI gap and flux recovery}\label{sec:recov}

For each problem separately, the constitutive residual is the correction
that a flux recovery makes. The flux $-(A\nabla w+\bff_2)$ of the
numerical potential is in general discontinuous, $\btau$ is a flux
computed in $H(\divv;\O)$, and $\Rcp(\btau,w)=\btau+A\nabla w+\bff_2$ is
the correction that $\btau$ makes to $-(A\nabla w+\bff_2)$, as a field,
before we take any norm.
Recovery-based estimators take a norm of it. They are the most widely
used a posteriori estimators in practice. The best known is the
Zienkiewicz--Zhu estimator \cite{ZZ:87,ZZ1:92}, which recovers a
continuous gradient by local averaging or patch projection. For
discontinuous coefficients the reconstruction should be
$H(\divv)$-conforming rather than continuous, since only the normal
component of the flux is continuous across interfaces \cite{CZ:09,CZ:10a}.
The least-squares recovery of \cite{CZ:10b} gives such a flux, and the
non-intrusive estimator of \cite{LZYZ:24} uses it. In this flux recovery
the recovered flux solves a problem posed in another space, $H(\divv)$.
It need not approximate $\bsigma$ better than the flux
$-(A\nabla w+\bff_2)$ of the potential. The reason is that the
least-squares functional, which adds the balance residual to
$\|A^{-1/2}\Rcp(\btau,w)\|_0^2$, is reliable and efficient by
\eqref{normeq} for arbitrary arguments, so we need no superconvergence
assumption.

The QoI gap uses the corrections themselves rather than their norms. By
\eqref{gap}, it pairs the primal and the adjoint flux corrections with
their signs, so it is a recovery-based estimate for the quantity of
interest. It has the same properties as this flux recovery: the
fluxes come from another space and need not be more accurate than the
gradients of the potentials, and \eqref{gap} and \eqref{csbound} hold for
arbitrary arguments. To combine two recovery estimators into a
goal-oriented indicator, a common choice is the sum over the elements of
the products of the local primal and adjoint indicators. This sum is
never negative, and by the Cauchy--Schwarz inequality it does not exceed
$\eta_p^{\LS}(\btau,w)\,\eta_d^{\LS}(\brho,v)$, which bounds the errors of the
evaluators up to a constant (Corollary~\ref{cor:bounds}). Unlike this
sum, the QoI gap keeps the signs of the two recovery fields.

\subsection{Symmetric problems and the primal--dual gap}\label{sec:pd}

For $\bb=\bzero$ the primal and adjoint problems are convex minimization
problems, and the primal--dual theory of \cite{Zhang:20} gives estimators
without any least-squares formulation. Let
$c\ge c_0>0$. The primal problem minimizes
$$
J_p(w):=\tfrac12(A\nabla w,\nabla w)+\tfrac12(cw,w)-(f_1,w)
+(\bff_2,\nabla w)+\langle\psi_p,w\rangle_{\G_N}
$$
over $X_p$, and its dual problem in the sense of convex duality, posed
over the fluxes, maximizes
$$
J_p^*(\btau):=-\tfrac12\big(A^{-1}(\btau+\bff_2),\btau+\bff_2\big)
-\tfrac12\big(c^{-1}(f_1-\divv\btau),f_1-\divv\btau\big)
-\langle\btau\cdot\bn,\phi_p\rangle_{\G_D}
$$
over all of $\bSigma_p$. We impose no equilibration constraint on
$\btau$. Integration by parts gives, for all $w\in X_p$ and
$\btau\in\bSigma_p$,
$$
J_p(w)-J_p^*(\btau)=\tfrac12\,\eta_{p,c}^{\LS}(\btau,w)^2
=\tfrac12\big(\|A^{1/2}\nabla e_p\|_0^2+\|c^{1/2}e_p\|_0^2\big)
+\tfrac12\big(\|A^{-1/2}E_p\|_0^2+\|c^{-1/2}\divv E_p\|_0^2\big),
$$
with $\eta_{p,c}^{\LS}$ of \eqref{LSc}. This is the generalized
Prager--Synge identity of \cite{Zhang:20}, and the same identity holds
for the adjoint problem. Since $a$ is then the energy inner product, the Cauchy--Schwarz
inequality, first in $L^2(\O)$ and then in $\mathbb R^2$, applied to
\eqref{gap_gamma0} and \eqref{err_gamma2}, gives the bound
$$
|\gamma-\gamma_0(w,v)|+|\gamma-\gamma_2(\btau,\brho)|\le\eta_{p,c}^{\LS}(\btau,w)\,\eta_{d,c}^{\LS}(\brho,v),
$$
and the bound of Proposition~\ref{prop:delta2} for $|\delta_2(w,\btau,v,\brho)|$ follows
from it. For $c=0$ the dual variables must satisfy
$\divv\btau=f_1$ and $\divv\brho=g_1$, that is, they must be
equilibrated, as in the bounds for outputs of \cite{PPP:97,SBHP:04,PNDP:21}.
Then $J_p(w)-J_p^*(\btau)=\frac12\|A^{-1/2}\Rcp(\btau,w)\|_0^2$, and in the same way
$$
|\gamma-\gamma_0(w,v)|+|\gamma-\gamma_1(w,\btau,v,\brho)|
\le\|A^{-1/2}\Rcp(\btau,w)\|_0\,\|A^{-1/2}\Rcd(\brho,v)\|_0.
$$
In these cases we do not need the least-squares functional. We can solve
the primal problem and its dual problem separately, by any methods,
while the least-squares functional couples the potential
and the flux. So we can take the flux from the dual problem, for
instance from a mixed method, which for $c=0$ gives the equilibrated flux
required there. If
we do not want to solve the mixed problem, we can use the least-squares
recovery instead, with the least-squares functional in place of the
primal--dual gap. For $c>0$, by
the identity above, minimizing $\eta_{p,c}^{\LS}(\btau,w)^2$ over $\btau\in\bSigma_p$ for a
fixed $w\in X_p$ maximizes $J_p^*$, so the weighted least-squares recovery
is the discrete dual problem and does not depend on $w$. We need the
least-squares functional when the problem is not a convex minimization
problem, as for $\bb\ne\bzero$, where no primal--dual gap is available.
For this reason we use the least-squares estimators in this paper. The QoI gap identity
\eqref{gap} and the bound \eqref{bound_gamma1} do not use the
minimization structure.

\section{Implementation}\label{sec:routes}

Theorem~\ref{thm:gap} makes no assumption on how we compute $(\btau,w)$
and $(\brho,v)$ or on how they are related. We can compute the QoI gap
whenever the four fields exist. For the bound we only evaluate the
least-squares functionals at these fields, so we do not need a
least-squares discretization for either. For the upper
bound to be good, however, all four fields must be good approximations.
By \eqref{normeq}, $\eta_p^{\LS}(\btau,w)$ is equivalent to the error of
the flux and the potential together, so a poor flux makes it large even
when the potential is accurate, and likewise for the adjoint pair. The
fluxes and the potentials should therefore have the same approximation
order. We obtain the fields from a least-squares finite element method,
which computes them directly, or, for an arbitrary discretization,
from the non-intrusive recovery of \cite{CZ:10b,LZYZ:24} or from another
method that computes a flux, such as a mixed method. The four fields must
also converge. The bounds hold for arbitrary arguments. They bound the
errors of the fields, but they do not show that these errors tend to
zero.
For a Galerkin potential with the recovered flux, refined with its own
estimator, \cite[\S7]{LZYZ:24} proves plain convergence of the adaptive
algorithm. The proof uses the a priori quasi-optimality of the Galerkin
method and of the recovery together with the reliability and efficiency
of the estimator. The a posteriori bounds alone do not give convergence.

Let $\cT_h$ and $\cT_H$ be conforming, shape-regular triangulations of
$\O$ into simplices, with elements $K$, for the primal and the adjoint
problem. The two triangulations need not be related. Throughout, $X_h\subset H^1(\O)$ denotes a $C^0$ Lagrange
space and $\bSigma_h\subset H(\divv;\O)$ a Raviart--Thomas space on
$\cT_h$, and $X_H$ and $\bSigma_H$ the corresponding spaces on $\cT_H$.
We write $S_k$ for the continuous piecewise polynomials of
degree $k$, $P_0$ for the piecewise constants and $RT_0$ for the
lowest-order Raviart--Thomas space. The affine subsets $X_h^p$ and
$\bSigma_h^p$ carry the discrete primal boundary conditions, $X_H^d$ and
$\bSigma_H^d$ the adjoint ones, and $\bSigma_h^0:=\bSigma_h\cap
H_N(\divv;\O)$ and $\bSigma_H^0:=\bSigma_H\cap H_N(\divv;\O)$ are the
linear subspaces.

\subsection{Fluxes from a least-squares finite element method}\label{sec:routeA}

If the underlying method is a first-order system least-squares method, the
fluxes are part of the discrete solution. We compute
\beq \label{lsfem}
(\bsigma_h,u_h)=\!\!\operatorname*{arg\,min}_{(\btau,w)\in
\bSigma_h^p\times X_h^p}\!\!\LS_p(\btau,w),
\qquad
(\br_H,z_H)=\!\!\operatorname*{arg\,min}_{(\brho,v)\in
\bSigma_H^d\times X_H^d}\!\!\LS_d(\brho,v),
\eeq
and insert them in \eqref{delta}. This is the setting of \cite{WZ:26},
which introduced $\gamma_1$, used it together with $\gamma_0$, and
proved bounds of the form \eqref{bound_gamma1}. The estimators $\eta_p^{\LS}(\bsigma_h,u_h)$ and
$\eta_d^{\LS}(\br_H,z_H)$ needed in Corollary~\ref{cor:cs} are the square
roots of the minimum values. The
same holds for any method whose unknowns are a conforming flux--potential
pair, such as the augmented mixed methods of \cite{MH:02,LZ:24}, with the
least-squares functionals evaluated at the computed pairs rather than
minimized.

\subsection{Non-intrusive least-squares flux recovery for an arbitrary potential}\label{sec:routeB}

A conforming Galerkin method, or any other method that returns only a
potential, gives no flux. We obtain one by minimizing the least-squares
functional over the flux alone, with the computed potential fixed: given $u_h\in X_h^p$ and $z_H\in X_H^d$, we set
\beq \label{recovery}
\bsigma_h:=\operatorname*{arg\,min}_{\btau\in\bSigma_h^p}
\LS_p(\btau,u_h),
\qquad
\br_H:=\operatorname*{arg\,min}_{\brho\in\bSigma_H^d}
\LS_d(\brho,z_H).
\eeq
This is the least-squares flux recovery \cite{CZ:10b,LZYZ:24} discussed
in Section~\ref{sec:recov}. The first-order conditions are
$(A^{-1}\Rcp(\bsigma_h,u_h),\bchi)+(\Rbp(\bsigma_h,u_h),\divv\bchi)=0$ for all $\bchi\in\bSigma_h^0$,
that is
\beq \label{recovery_ELp}
b_{\rm flux}(\bsigma_h,\bchi)
= -\big(A^{-1}(\bff_2+A\nabla u_h),\,\bchi\big)
- \big(\bb\cdot\nabla u_h+c\,u_h-f_1,\;\divv\bchi\big),
\eeq
and, for all $\bchi\in\bSigma_H^0$,
\beq \label{recovery_ELd}
b_{\rm flux}(\br_H,\bchi)
= -\big(A^{-1}(\bg_2+A\nabla z_H+\bb\,z_H),\,\bchi\big)
- \big(c\,z_H-g_1,\;\divv\bchi\big),
\eeq
with, in both cases, the $A$-weighted $H(\divv)$ inner product
\beq \label{bflux}
b_{\rm flux}(\bmu,\bchi):=\big(A^{-1}\bmu,\bchi\big)
+\big(\divv\bmu,\divv\bchi\big),
\qquad b_{\rm flux}(\bmu,\bmu)=\tri(\bmu,0)\tri^2 .
\eeq
Both problems are therefore uniquely solvable. The potentials enter only
the right-hand sides, so we can use potentials in $H^1(\O)$ from any
method and on any mesh in place of $u_h$ and $z_H$. The error of the recovered flux in $H(\divv;\O)$ is
bounded by the best approximation error of $\bsigma$ in $\bSigma_h^p$
plus a constant times the $H^1$ error of $u_h$. This is \cite[Thm.~5.1]{LZYZ:24} for
$\G_N=\emptyset$, and the following argument covers $\G_N\ne\emptyset$.
Since $\Rcp(\bsigma,u)=\Rbp(\bsigma,u)=0$, the exact flux satisfies
\eqref{recovery_ELp} with $u_h$ replaced by $u$, so that
$$
b_{\rm flux}(\bsigma-\bsigma_h,\bchi)=-\big(\nabla(u-u_h),\bchi\big)
-\big(\bb\cdot\nabla(u-u_h)+c\,(u-u_h),\,\divv\bchi\big)
\qquad\forall\,\bchi\in\bSigma_h^0 .
$$
Let $P_h\bsigma\in\bSigma_h^p$ be the $b_{\rm flux}$-orthogonal projection
of $\bsigma$, that is, $b_{\rm flux}(\bsigma-P_h\bsigma,\bchi)=0$ for all
$\bchi\in\bSigma_h^0$. We test with $\bchi=P_h\bsigma-\bsigma_h$ and obtain
$\tri(P_h\bsigma-\bsigma_h,0)\tri\le C\,\tri(0,u-u_h)\tri$, with $C$
depending on $A$, $\bb$ and $c$. Since $\tri(\bsigma-P_h\bsigma,0)\tri$ is
the best approximation error, the triangle inequality gives the bound.
This bound is of the lower of the two approximation orders, which is
another reason to choose the flux space and the potential space of the
same approximation order. We use $RT_0$ with $S_1$ throughout.

The flux need not come from the least-squares recovery. We can use any flux in
$\bSigma_p$, for instance that of a mixed method for the
primal problem, and likewise for the adjoint problem. For $\bb=\bzero$
this is the setting of Section~\ref{sec:pd}, in which the primal problem
and its dual problem are decoupled and we compute the potential and the
flux separately.

\begin{rem}[Recovering the potential]\label{rem:potrec}
For a method that delivers a flux but no conforming potential, such as a
mixed method, we can obtain a potential in the same way, by minimizing
$\LS_p(\bsigma_h,\cdot)$ over $X_h^p$. This is a second-order problem on
$X_h$ with the symmetric positive definite form
$(A\nabla w,\nabla s)+(\bb\cdot\nabla w+cw,\bb\cdot\nabla s+cs)$. There are local alternatives. A direct interpolation of the
discontinuous potential gives a conforming one. Where only the gradient
is needed, as in the constitutive residual $\Rcp(\btau,w)$, a recovery in
$H(\operatorname{curl})$ as in \cite{CCZ:20mixed} gives an approximation
of $\nabla u$ without a potential, but the evaluators also need $w$
itself.
\end{rem}

\subsection{The adaptive procedure}\label{sec:procedure}

On each mesh we compute the primal and adjoint approximations
$(\bsigma_h,u_h)$ and $(\br_h,z_h)$, by a least-squares finite element
method or by another method followed by the flux recovery
\eqref{recovery}, together with the least-squares estimators. When we
solve both problems on the same mesh, we mark elements by D\"orfler's
strategy applied to the balanced indicator of \cite{WZ:26},
\beq \label{balanced}
\eta_K^2:=\frac12\Big(\frac{\eta_d^{\LS}}{\eta_p^{\LS}}\big(\eta_{p,K}^{\LS}\big)^2
+\frac{\eta_p^{\LS}}{\eta_d^{\LS}}\big(\eta_{d,K}^{\LS}\big)^2\Big),
\qquad K\in\cT_h,
\eeq
where $\eta_p^{\LS}=\eta_p^{\LS}(\bsigma_h,u_h)$,
$\eta_d^{\LS}=\eta_d^{\LS}(\br_h,z_h)$, and $\eta_{p,K}^{\LS}$ and
$\eta_{d,K}^{\LS}$ are the square roots of the least-squares functionals
\eqref{LS} restricted to $K$. For recovered pairs these are the local
non-intrusive estimators of \cite{LZYZ:24}. The weights in
\eqref{balanced} are the optimal ones in Young's inequality applied to the
product $\eta_p^{\LS}\eta_d^{\LS}$, so that
$\sum_{K\in\cT_h}\eta_K^2=\eta_p^{\LS}\eta_d^{\LS}$, similar to the weighted
marking of \cite{BET:11}. On each element the balanced indicator is at
least the product of the local indicators, and it is positive on every
element where either of them is positive.

We stop when the difference of the two evaluators that we report is below
a tolerance $\mathrm{TOL}$ chosen by the user. We report the interval
between the two evaluators and its midpoint rather than either evaluator,
because the computed quantities do not tell which of the two is the more
accurate. For $\gamma_0$ and $\gamma_1$ the criterion is
$|\delta(u_h,\bsigma_h,z_h,\br_h)|\le\mathrm{TOL}$. Since
$|\delta(u_h,\bsigma_h,z_h,\br_h)|\le\eta_p^{\LS}(\bsigma_h,u_h)\,\eta_d^{\LS}(\br_h,z_h)$
by Corollary~\ref{cor:cs}, the criterion holds after finitely many steps
of any refinement under which this product tends to zero. Where reaction dominates, the signed QoI gap does not contain the
reaction term (Section~\ref{sec:gamma2}) and can be small while the errors of
$\gamma_0(u_h,z_h)$ and $\gamma_1(u_h,\bsigma_h,z_h,\br_h)$ are large and
nearly equal. There we report $\gamma_0(u_h,z_h)$ and
$\gamma_2(\bsigma_h,\br_h)$ and use the criterion
$|\delta_2(u_h,\bsigma_h,z_h,\br_h)|\le\mathrm{TOL}$,
which \eqref{delta2bound} makes reachable in the same way. We do not use
$\gamma_2$ where $\vartheta$ of Proposition~\ref{prop:delta2} is large,
since there the terms weighted by $c^{-1}$ can make
$\gamma_2(\bsigma_h,\br_h)$ less accurate than
$\gamma_1(u_h,\bsigma_h,z_h,\br_h)$ (Remark~\ref{rem:gamma2}). Proposition~\ref{prop:twosided} relates the
length of the reported interval to the errors of its endpoints.

\section{Numerical experiments}\label{sec:num}

\subsection{Setup}\label{sec:setup}

We compute in two dimensions, with lowest-order elements on
triangulations that we refine adaptively by newest-vertex bisection. Only
the comparison with finite differences uses uniform grids. We use two kinds of
pairs $(\bsigma_h,u_h)$ and $(\br_h,z_h)$.
\begin{itemize}
\item[LS:] the least-squares method \eqref{lsfem} with $RT_0\times S_1$,
for the primal and the adjoint problem (Section~\ref{sec:routeA});
\item[G:] conforming $S_1$ Galerkin approximations of $u$ and $z$, with
both fluxes in $RT_0$ from the least-squares recovery \eqref{recovery}
(Section~\ref{sec:routeB}).
\end{itemize}
We also combine the G potentials with the fluxes of the mixed method
$RT_0\times P_0$, and use pairs from finite difference and finite volume
schemes and pairs from different methods or meshes. On each mesh we write
$$
\gamma_{0,h}:=\gamma_0(u_h,z_h),\qquad
\gamma_{1,h}:=\gamma_1(u_h,\bsigma_h,z_h,\br_h),\qquad
\delta_h:=\delta(u_h,\bsigma_h,z_h,\br_h)=\gamma_{1,h}-\gamma_{0,h},
$$
and $\eta_p^{\LS}$, $\eta_d^{\LS}$, $\eta_{p,c}^{\LS}$, $\eta_{d,c}^{\LS}$
for the estimators at $(\bsigma_h,u_h)$ and $(\br_h,z_h)$. Unless stated
otherwise, we solve the primal and adjoint problems on the same mesh
sequence and refine by the procedure of Section~\ref{sec:procedure}, with
D\"orfler bulk parameter $0.5$. We compute
$\delta_h$ by subtraction or by quadrature of the right-hand side of
\eqref{delta}. With the same quadrature rule for corresponding terms, the
two agree up to rounding. Table~\ref{tab:ref} lists the reference values $\gamma$ and how we
compute them. On each edge of $\G_N$ we set the normal component of the
fluxes in $RT_0$ to the value of $\psi_p$ or $\psi_d$ at the midpoint of
the edge. We interpolate the
Dirichlet data $\phi_p$ and $\phi_d$ at the boundary nodes, so in general
$u_h\notin X_p$ and $z_h\notin X_d$.
For $c>0$ we let $\gamma_{2,h}$ be the last expression in
\eqref{gamma2res} at the computed fields and
$\delta_{2,h}:=\gamma_{2,h}-\gamma_{0,h}$. For $u_h\in X_p$ and
$z_h\in X_d$ these are $\gamma_2(\bsigma_h,\br_h)$ and
$\delta_2(u_h,\bsigma_h,z_h,\br_h)$. The errors $\gamma-\gamma_{i,h}$,
$i=0,1,2$, then contain the same boundary term $\langle\br_h\cdot\bn,\phi_p-u_h\rangle_{\G_D}
+\langle\bsigma_h\cdot\bn,\phi_d-z_h\rangle_{\G_D}$, which cancels in the
QoI gaps $\delta_h$ and $\delta_{2,h}$. We leave its treatment to the
next paper of this series. We measure complexity by the number $N$ of vertices of the mesh.
The tolerance $\mathrm{TOL}$ in the stopping criterion of
Section~\ref{sec:procedure} is the user's choice, and the criterion needs
$|\delta_h|$ to be of the size of the larger of the errors of
$\gamma_{0,h}$ and $\gamma_{1,h}$ on whatever mesh the computation
stops. For this reason we do not stop the computations at a particular
tolerance. We run them to a prescribed size, given with each experiment,
and compare $|\delta_h|$ with the errors on every mesh.

We report the signed errors $\gamma-\gamma_{0,h}$ and
$\gamma-\gamma_{1,h}$ of the two evaluators, the QoI gap $\delta_h$, the
midpoint $\bar\gamma_h:=\frac12(\gamma_{0,h}+\gamma_{1,h})$, the product
bound $\eta_p^{\LS}\eta_d^{\LS}$ and the effectivity index
\beq \label{effs}
\mathrm{eff}:=\frac{|\delta_h|}{\max(|\gamma-\gamma_{0,h}|,|\gamma-\gamma_{1,h}|)} ,
\eeq
which is the quantity controlled by Proposition~\ref{prop:twosided}.
To compare with it, we report the factor
$\eta_p^{\LS}\eta_d^{\LS}/|\gamma-\gamma_{0,h}|$ by which the
product bound exceeds the error of $\gamma_{0,h}$. We also record, a posteriori, which case of
Proposition~\ref{prop:twosided} occurs, and the ratio of the smaller to the larger error,
$$
\mathrm{ratio}:=\frac{\min(|\gamma-\gamma_{0,h}|,|\gamma-\gamma_{1,h}|)}{\max(|\gamma-\gamma_{0,h}|,|\gamma-\gamma_{1,h}|)},
$$
which is $|e_A|/|e_B|$ in Proposition~\ref{prop:twosided}. Then
$\mathrm{eff}=1+\mathrm{ratio}$ in case (a) and $\mathrm{eff}=1-\mathrm{ratio}$
in case (b). We take all entries of a table row from the same mesh, so
that every row satisfies $\delta_h=\gamma_{1,h}-\gamma_{0,h}$.

\subsection{Test problems P1--P7}\label{sec:tests}

We take problems P1--P6 from \cite{WZ:26}, which shows the meshes and the
convergence histories of the adaptive computations. We add P7 to compare
equilibrated and recovered fluxes. We compute
the data not stated from the exact solutions: $f_1$, $\phi_p$ and
$\psi_p$ from $u$ by \eqref{primal_pde}, and, in P2, $g_1$,
$\phi_d$ and $\psi_d$ from $z$ by \eqref{dual_pde}. All other data are
zero.

\begin{enumerate}
\item[P1.] {\em Mommer--Stevenson benchmark} \cite{MS:09}.
$\O=(-1,1)^2$, $\G_D=\partial\O$, $A=I$, $\bb=\bzero$, $c=0$,
$\bff_2=(\mathbf 1_{T_p},0)^T$ and $\bg_2=(\mathbf 1_{T_d},0)^T$, where
$T_p$ is the triangle with vertices $(1,0),(1,1),(0,1)$ and $T_d$ the
triangle with vertices $(-1,-1),(0,-1),(-1,0)$. The primal and adjoint
residuals concentrate near the two disjoint triangles, so a local
product indicator is nearly zero everywhere. We do not know the exact
solution in closed form, and Table~\ref{tab:ref} describes how we compute
the reference value.

\item[P2.] {\em Boundary layers.} $\O=(0,1)^2$, $\G_D=\partial\O$,
$A=\epsilon I$, $\bb=(b_1,b_2)^T=(2,3)^T$, $c=1$, and
$$
u=\big(x-e^{b_1(x-1)/\epsilon}\big)\big(y^2-e^{b_2(y-1)/\epsilon}\big),
\qquad
z=\big(1-x-e^{-b_1x/\epsilon}\big)\big((1-y)^2-e^{-b_2y/\epsilon}\big),
$$
with layers at $x=1$ and $y=1$ in $u$ and at $x=0$ and $y=0$ in $z$. Here
$\epsilon=0.1$, except in the computations with varying $c$, where
$\epsilon=0.05$. The adjoint Dirichlet datum $\phi_d=z|_{\partial\O}$ is of size
$e^{-b_1/\epsilon}$.

\item[P3--P6.] {\em L-shaped domain.}
$\O=(-1,1)^2\setminus([0,1]\times[-1,0])$, $\G_N=\{(x,1):-1<x<1\}$,
$\G_D=\partial\O\setminus\overline\G_N$, $A=I$, $\bb=\bzero$, $c=0$,
and, in polar coordinates $(r,\theta)$ about the origin,
$u=r^{2/3}\sin(2\theta/3)$, so that $\bsigma=-\nabla u$. The four
problems share the primal problem and differ in the adjoint data: in P3,
$g_1=\mathbf 1_{\omega_0}$ with $\omega_0=(-1,-\frac12)^2$; in P4,
$\bg_2=(\mathbf 1_{T_0},0)^T$ with $T_0$ the triangle with vertices
$(-1,-1),(-1,0),(-\frac12,-\frac12)$; in P5, $\psi_d=-1$; in P6,
$\phi_d=y+1$ on $\G_L=\{(-1,y):-1<y<1\}$ and $\phi_d=0$ on
$\G_D\setminus\G_L$.

\item[P7.] {\em Kellogg's problem} \cite{Kellogg:74}, as in \cite{CZ:09}.
$\O=(-1,1)^2$, $\G_D=\partial\O$, $\bb=\bzero$, $c=0$, $A=aI$ with
$a=R\approx161.4476387975881$ in the first and third quadrants and $a=1$
in the second and fourth, and $u=r^{\beta}\mu(\theta)$ with $\beta=0.1$
and $\mu$ as in \cite{CZ:09}. The solution satisfies $u\in H^{1+t}(\O)$
only for $t<\beta$.
The adjoint data are $\bg_2=(\mathbf 1_{\omega},0)^T$ with
$\omega=(-\frac12,\frac12)\times(0,\frac12)$, which straddles the
interface $x=0$, so that
$\gamma=-(\partial_xu,\mathbf 1_{\omega})=(A^{-1}\bsigma,\bg_2)$.
\end{enumerate}

\begin{table}[htbp]
\small
\begin{tabular}{llr}
\toprule
& quantity of interest & reference value $\gamma$ \\
\midrule
P1 & $-(\partial_xu,\mathbf 1_{T_d})$
   & $-6.3403633\cdot10^{-3}$ \\
P2 & $(g_1,u)+\langle\bsigma\cdot\bn,\phi_d\rangle_{\partial\O}$
   & $1.0528716\cdot10^{-1}$ \\
P3 & $(u,\mathbf 1_{\omega_0})$ & $1.3000215\cdot10^{-1}$ \\
P4 & $(\sigma_1,\mathbf 1_{T_0})$ & $1.5771219\cdot10^{-1}$ \\
P5 & $\langle u,1\rangle_{\G_N}$ & $1.7886031$ \\
P6 & $\langle\bsigma\cdot\bn,\,y+1\rangle_{\G_L}$ & $-1.0326505$ \\
P7 & $(A^{-1}\bsigma,\bg_2)$
   & $5.3272433\cdot10^{-2}$ \\
\bottomrule
\end{tabular}
\caption{Quantities of interest and reference values. For P2--P7
we compute them by quadrature of the exact solutions on the subdomains
where these are smooth. For P1 we use $\gamma_0(w,v)$ with
$S_3$ and with $S_4$ Galerkin approximations of the
primal and adjoint problems, on meshes graded towards the four points
where the hypotenuses of $T_p$ and $T_d$ meet $\partial\O$. The two
values agree to $4\cdot10^{-14}$.}
\label{tab:ref}
\end{table}

\subsection{Least-squares pairs}\label{sec:num_main}

\begin{figure}[htbp]
\centering
\includegraphics[width=.6\textwidth]{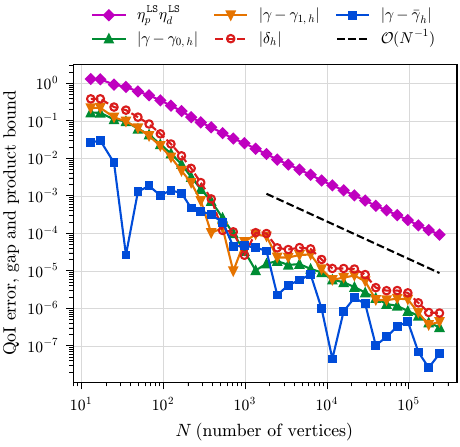}
\caption{LS pairs for P2: the errors $|\gamma-\gamma_{0,h}|$, $|\gamma-\gamma_{1,h}|$ of the two evaluators and
$|\gamma-\bar\gamma_h|$ of their midpoint, the QoI gap $|\delta_h|$ and the
product bound $\eta_p^{\LS}\eta_d^{\LS}$ against the number $N$ of
vertices.}
\label{fig:conv}
\end{figure}

Table~\ref{tab:main} shows the values for LS pairs on the final mesh of
each sequence, and Figure~\ref{fig:conv} the whole sequence of P2. We
refine until the number of vertices first exceeds $2\cdot10^5$.

\begin{table}[htbp]
\small
\setlength{\tabcolsep}{4.5pt}
\begin{tabular}{lrccccccccc}
\toprule
& $N$ & $\gamma-\gamma_{0,h}$ & $\gamma-\gamma_{1,h}$ & $\delta_h$ & $\mathrm{ratio}$ & case & $\mathrm{eff}$
& $\dfrac{\eta_p^{\LS}\eta_d^{\LS}}{|\gamma-\gamma_{0,h}|}$
& $|\gamma-\bar\gamma_h|$ & $\eta_p^{\LS}\eta_d^{\LS}$ \\
\midrule
P1 & $205083$ & $-2.29$e-7 & $9.97$e-8 & $-3.29$e-7 & $0.44$ & (a) & $1.44$ & $18.2$ & $6.46$e-8 & $4.16$e-6 \\
P2 & $236439$ & $-3.14$e-7 & $4.37$e-7 & $-7.51$e-7 & $0.72$ & (a) & $1.72$ & $295$ & $6.17$e-8 & $9.26$e-5 \\
P3 & $241295$ & $2.74$e-8 & $-4.91$e-8 & $7.65$e-8 & $0.56$ & (a) & $1.56$ & $47.6$ & $1.09$e-8 & $1.30$e-6 \\
P4 & $251015$ & $-7.47$e-7 & $5.98$e-8 & $-8.07$e-7 & $0.08$ & (a) & $1.08$ & $5.98$ & $3.44$e-7 & $4.47$e-6 \\
P5 & $236900$ & $1.74$e-6 & $-1.14$e-6 & $2.88$e-6 & $0.65$ & (a) & $1.65$ & $7.33$ & $3.02$e-7 & $1.27$e-5 \\
P6 & $250369$ & $3.98$e-6 & $-2.28$e-6 & $6.26$e-6 & $0.57$ & (a) & $1.57$ & $3.25$ & $8.51$e-7 & $1.29$e-5 \\
\bottomrule
\end{tabular}
\caption{LS pairs, final mesh of each sequence. $N$ is the number of
vertices, and we compute all entries of a row on the same mesh.}
\label{tab:main}
\end{table}

On the final meshes case (a) occurs in every problem. Along the
sequences, case (a) occurs on every mesh of P1 and P3--P6,
with $\mathrm{eff}$ in $[1.23,1.97]$, $[1.29,1.96]$, $[1.03,1.49]$,
$[1.51,1.98]$ and $[1.45,1.75]$, respectively. In P2, case (a) occurs on
$31$ of the $33$ meshes, with $\mathrm{eff}\in[1.10,2.00]$, and case (b)
on meshes $13$ and $15$, with $\mathrm{eff}=0.46$ and $0.44$. From mesh
$12$ to mesh $15$, $\gamma-\gamma_{1,h}$ changes sign three times, and
$\gamma-\gamma_{0,h}$ then changes sign between meshes $15$ and $16$. The
QoI gap changes sign once, between meshes $14$ and $15$, and in the other
problems it does not change sign.

Over the last six meshes, the fitted rates against $N$ in P1--P6 lie
between $-0.50$ and $-0.51$ for $\eta_p^{\LS}$ and $\eta_d^{\LS}$, between
$-1.00$ and $-1.01$ for $\eta_p^{\LS}\eta_d^{\LS}$, and between $-0.99$
and $-1.11$ for $|\delta_h|$. For $|\gamma-\gamma_{0,h}|$ and
$|\gamma-\gamma_{1,h}|$ they lie between $-0.85$ and $-1.18$. On the
final meshes $\mathrm{eff}$ lies between $1.08$ and $1.72$, while the
product $\eta_p^{\LS}\eta_d^{\LS}$ exceeds $|\gamma-\gamma_{0,h}|$ by
factors between $3.2$ and $295$ (Table~\ref{tab:main}). In P2 this factor
lies between $258$ and $315$ on the last ten meshes. The
product exceeds $|\delta_h|$ by
the factors $13$ and $123$ in P1 and P2 and $17$, $5.5$, $4.4$ and
$2.1$ in P3--P6. On the last ten meshes this factor stays within $10\%$
of its final value in P1 and P3--P6 and lies between $88$ and $157$ in
P2. The product and $|\delta_h|$ both decrease like $N^{-1}$. If we
stopped when the product, instead of $|\delta_h|$, fell below a
tolerance, we would need about these factors more vertices, $123$ times
as many in P2. Of P1--P6, only P2 is convection-dominated, and it is the
only problem in which $\eta_p^{\LS}\eta_d^{\LS}$ exceeds
$|\gamma-\gamma_{0,h}|$ by more than a factor $48$.

P2, where $c=1$ and $\vartheta=11.4$, shows why we do not use $\gamma_2$
when $\vartheta$ is large. There $\gamma_{2,h}$ is the least accurate
evaluator on every mesh, by the factors $55$ and $40$ against
$\gamma_{0,h}$ and $\gamma_{1,h}$ on the final mesh. The pair $\gamma_{0,h},\gamma_{2,h}$ is
in case (b) on the first $15$ meshes and in case (a) on the remaining
$18$. On these $18$ meshes $\mathrm{ratio}\le0.023$ and
$\mathrm{eff}\in[1.003,1.023]$, so that $|\delta_{2,h}|$ is the error
of $\gamma_{2,h}$ to within $2.3\%$. On the final mesh
$\gamma-\gamma_{2,h}=1.73\cdot10^{-5}$ and $\delta_{2,h}=-1.76\cdot10^{-5}$, and the product
$\eta_{p,c}^{\LS}\eta_{d,c}^{\LS}$, here equal to
$\eta_p^{\LS}\eta_d^{\LS}$, exceeds $|\delta_{2,h}|$ by the factor $5.3$,
against $123$ for $|\delta_h|$.

\subsection{Other methods, other fluxes and unrelated pairs}\label{sec:num_rec}

In this subsection at least one of the two pairs is not a least-squares
pair on the common mesh. Either the potentials come from a Galerkin,
finite difference or finite volume method and we recover the fluxes or
take them from a mixed method, or the two pairs come from different
methods or meshes.
Theorem~\ref{thm:gap} and Proposition~\ref{prop:twosided} apply
unchanged. We test whether the QoI gap stays accurate, and how this depends on
the method that produced the potentials, on how we obtain the fluxes,
and on the relation between the two pairs.

\subsubsection{Conforming Galerkin potentials}\label{sec:num_G}

For G pairs on a common mesh and $\phi_d=0$, $\gamma_{0,h}=m(u_h)$
(Section~\ref{sec:dwr}) is the value of the quantity of interest computed
directly from the Galerkin solution, and the recovered fluxes enter only
through $\delta_h$, since $\gamma_{1,h}=m(u_h)+\delta_h$. In P6, where
$m=0$ and $\phi_d\ne0$, $\gamma_{0,h}=\ell(z_h)-a(u_h,z_h)$ depends only on
the boundary values of $z_h$ and is the variationally consistent
evaluation of the boundary flux. A dual-weighted residual estimate of
$\gamma-\gamma_{0,h}$ needs a second, enriched adjoint approximation, since by Galerkin
orthogonality \eqref{etadwr} vanishes for every $\tilde v$ in the
Galerkin space with the boundary values of $z_h$. The QoI gap needs the
two recovered fluxes instead.

\begin{table}[htbp]
\small
\setlength{\tabcolsep}{4.5pt}
\begin{tabular}{lrccccccccc}
\toprule
& $N$ & $\gamma-\gamma_{0,h}$ & $\gamma-\gamma_{1,h}$ & $\delta_h$ & $\mathrm{ratio}$ & case & $\mathrm{eff}$
& $\dfrac{\eta_p^{\LS}\eta_d^{\LS}}{|\gamma-\gamma_{0,h}|}$
& $|\gamma-\bar\gamma_h|$ & $\eta_p^{\LS}\eta_d^{\LS}$ \\
\midrule
P1 & $274402$ & $-1.31$e-7 & $1.24$e-7 & $-2.56$e-7 & $0.95$ & (a) & $1.95$ & $24.1$ & $3.46$e-9 & $3.17$e-6 \\
P2 & $234550$ & $-3.14$e-7 & $4.37$e-7 & $-7.51$e-7 & $0.72$ & (a) & $1.72$ & $297$ & $6.13$e-8 & $9.34$e-5 \\
P3 & $230503$ & $2.89$e-8 & $-5.13$e-8 & $8.02$e-8 & $0.56$ & (a) & $1.56$ & $47.1$ & $1.12$e-8 & $1.36$e-6 \\
P4 & $251207$ & $-7.47$e-7 & $5.96$e-8 & $-8.06$e-7 & $0.08$ & (a) & $1.08$ & $5.98$ & $3.44$e-7 & $4.46$e-6 \\
P5 & $232596$ & $1.78$e-6 & $-1.18$e-6 & $2.96$e-6 & $0.66$ & (a) & $1.66$ & $7.29$ & $3.04$e-7 & $1.30$e-5 \\
P6 & $240333$ & $4.22$e-6 & $-2.19$e-6 & $6.41$e-6 & $0.52$ & (a) & $1.52$ & $3.20$ & $1.02$e-6 & $1.35$e-5 \\
\bottomrule
\end{tabular}
\caption{G pairs with least-squares flux recovery, final mesh of each
sequence. We generate the meshes as for LS pairs, with the indicator
\eqref{balanced} evaluated at the G pairs.}
\label{tab:G}
\end{table}

We compare the G pairs of Table~\ref{tab:G} with the LS pairs of
Table~\ref{tab:main}. On the final meshes the case is the same in every
problem, and the factor by which the product exceeds the QoI gap agrees to
within $3\%$. The ratio of the two errors agrees to within $0.06$
in all problems except P1, where it is not stable along either sequence.
On the last $10$ meshes of P1, $|\gamma-\gamma_{0,h}|/|\gamma-\gamma_{1,h}|$
ranges over $[0.41,2.4]$ for G pairs and over $[0.41,2.3]$ for LS pairs. Along the sequences, case (a) occurs on every
mesh of P1, P2, P3, P5 and P6, and on $30$ of the $31$ meshes of P4. The
exception is mesh $12$ of P4, at which the smaller error
$\gamma-\gamma_{1,h}$ changes sign, with $\mathrm{ratio}=0.10$ and
$\mathrm{eff}=0.90$. On this mesh $|\delta_h|$ underestimates $\max(|\gamma-\gamma_{0,h}|,|\gamma-\gamma_{1,h}|)$ by $10\%$. In P2 the
QoI gap of the G pairs keeps its sign along the whole sequence, and no mesh
of the G sequence is in case (b). As for LS pairs, $\gamma_{0,h}$ is the more accurate on most meshes of
P2, here on $27$ of the $34$.
With $\phi_d=0$, in P1 and P3--P5, the computed $\gamma_{0,h}$ agrees with
$m(u_h)$ on the final meshes to within $3\cdot10^{-6}|\gamma-\gamma_{0,h}|$.

\subsubsection{Equilibrated and recovered fluxes for Kellogg's problem}
\label{sec:num_kellogg}

For P7 we combine the G potentials $u_h$ and $z_h$ on each mesh with two
kinds of fluxes in $RT_0$:
\begin{itemize}
\item[(i)] the fluxes of the mixed method $RT_0\times P_0$ for the primal
and the adjoint problem, which are equilibrated, $\divv\bsigma_h=f_1$ and
$\divv\br_h=g_1$, both zero here;
\item[(ii)] the least-squares recovery \eqref{recovery} \cite{CZ:10b}.
\end{itemize}
We generate the meshes with the pairs (i) and the indicator
\eqref{balanced}. Both fluxes are $H(\divv)$-conforming, as discontinuous
coefficients require \cite{CZ:09,CZ:12}, so
Theorem~\ref{thm:gap} applies to each. On each mesh $\gamma_{0,h}$ is
the same for (i) and (ii), and the two QoI gaps differ only through the
fluxes. We take the LS pair on the same meshes as a third variant. The balance
residuals of the pairs (i) vanish here, and since $\bb=\bzero$ and $c=0$,
P7 with the fluxes
(i) is the equilibrated symmetric case of Section~\ref{sec:pd}.
We ask whether the accuracy of the QoI gap depends on how we obtain the
flux.

\begin{figure}[htbp]
\centering
\includegraphics[width=.5\textwidth]{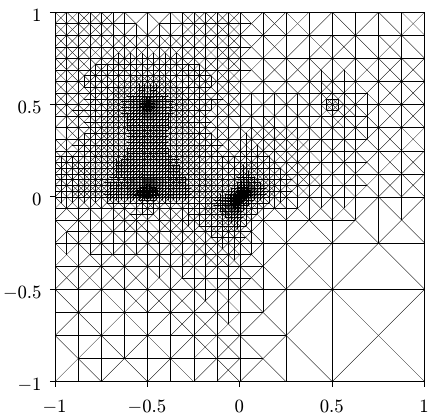}
\caption{P7: adaptive mesh $19$, with $N=3456$. The refinement
concentrates at the origin and near the corners of $\omega$.}
\label{fig:kellogg}
\end{figure}

We compute the sequence up to mesh $28$, with $N=22856$, and
Figure~\ref{fig:kellogg} shows mesh $19$. Meshes refined for the energy
error of $u$ are radial about the origin \cite{Zhang:26k}. Our meshes also
resolve the adjoint solution, because the indicator \eqref{balanced}
distributes $\eta_p^{\LS}\eta_d^{\LS}$ over the elements. The adjoint data
$\bg_2=(\mathbf 1_{\omega},0)^T$ are discontinuous across the vertical
edges of $\omega$, so $z$ is singular at the corners of $\omega$. Near the
corners in the second quadrant, where $a=1$, the meshes are much finer
than near those in the first, where $a=R$. The three
variants give nearly the same values on every mesh. Case (a) occurs on
every mesh, with $\mathrm{eff}\in[1.08,1.84]$. On mesh $28$ the three
rows of Table~\ref{tab:kellogg} agree in the digits shown, except
$\gamma-\gamma_{1,h}$ and $|\gamma-\bar\gamma_h|$ in the third digit, and
the exact errors of the primal fluxes agree as well,
$(\|A^{-1/2}E_p\|_0^2+\|\divv E_p\|_0^2)^{1/2}=0.171$. On this problem the
accuracy of the QoI gap does not depend on how we obtain the flux. The
product exceeds the QoI gap only by factors between $1.7$ and $2.7$ on meshes
$19$--$28$, as in P6. The
midpoint is more accurate than $\gamma_{0,h}$ by the factor $12$. The
errors have not yet reached their asymptotic rates on this sequence. Over
meshes $23$--$28$ the fitted rates
against $N$ are $-0.15$ for the flux error, $-0.58$ for
$\eta_p^{\LS}\eta_d^{\LS}$, and $-0.37$ and $-0.41$ for $|\delta_h|$ and
$\max(|\gamma-\gamma_{0,h}|,|\gamma-\gamma_{1,h}|)$.

\begin{table}[htbp]
\small
\setlength{\tabcolsep}{4.5pt}
\begin{tabular}{lccccccccc}
\toprule
& $\gamma-\gamma_{0,h}$ & $\gamma-\gamma_{1,h}$ & $\delta_h$ & $\mathrm{ratio}$ & case & $\mathrm{eff}$
& $\dfrac{\eta_p^{\LS}\eta_d^{\LS}}{|\gamma-\gamma_{0,h}|}$
& $|\gamma-\bar\gamma_h|$ & $\eta_p^{\LS}\eta_d^{\LS}$ \\
\midrule
G, (i) mixed flux    & $7.34$e-4 & $-6.17$e-4 & $1.35$e-3 & $0.84$ & (a) & $1.84$ & $3.13$ & $5.87$e-5 & $2.30$e-3 \\
G, (ii) LS recovery  & $7.34$e-4 & $-6.16$e-4 & $1.35$e-3 & $0.84$ & (a) & $1.84$ & $3.13$ & $5.89$e-5 & $2.30$e-3 \\
LS                   & $7.34$e-4 & $-6.16$e-4 & $1.35$e-3 & $0.84$ & (a) & $1.84$ & $3.13$ & $5.89$e-5 & $2.30$e-3 \\
\bottomrule
\end{tabular}
\caption{P7 on the mesh sequence generated with the pairs (i), mesh $28$,
$N=22856$. For G pairs, $\gamma-\gamma_{0,h}$ is the same in both rows.}
\label{tab:kellogg}
\end{table}

\subsubsection{Pairs from different methods and meshes, finite differences and finite volumes}\label{sec:num_indep}

Theorem~\ref{thm:gap} requires no relation between the two pairs
and nothing of the methods that produced them. We
compute the QoI gap from
\begin{itemize}
\item[G--LS:] a G primal pair and an LS adjoint pair on a common mesh, for
P1--P6, with the marking of Section~\ref{sec:setup};
\item[LS--LS:] LS primal and adjoint pairs on separate mesh sequences, each
refined by D\"orfler marking of its own least-squares functional, for P1
and P2. We pair the $k$-th meshes of the two sequences and evaluate
$\gamma_{0,h}$, $\gamma_{1,h}$ and $\delta_h$ on their overlay, the coarsest common refinement, on which both pairs are
piecewise polynomial;
\item[FD:] for P2, pairs from a finite difference scheme on uniform grids
of mesh size $h=2^{-3},\dots,2^{-9}$, without adaptivity. We apply the
five-point scheme with upwinding of the convection term to the primal
problem and to the adjoint problem \eqref{dual_pde}, in the latter with
upwinding for $-\bb$. We take the grid values as a continuous piecewise
linear function on the triangulation that cuts every grid square along
the same diagonal, and obtain both fluxes in $RT_0$ by the recovery
\eqref{recovery}, as for G pairs. For comparison we also compute G pairs
on the same grids;
\item[FV:] for P1--P7, potentials from the vertex-centered finite volume
method (box method) on the barycentric dual mesh, with upwinding of the
convection term, for the primal problem and for the adjoint problem
\eqref{dual_pde} in conservative form. We take the nodal values as
functions in $S_1$, obtain both fluxes by the recovery \eqref{recovery},
and generate the meshes as for G pairs, with the indicator
\eqref{balanced} evaluated at the finite volume pairs.
\end{itemize}
Table~\ref{tab:indep} shows the results on the final meshes.

\begin{table}[htbp]
\small
\setlength{\tabcolsep}{4.5pt}
\begin{tabular}{llrccccccccc}
\toprule
& & $N$ & $\gamma-\gamma_{0,h}$ & $\gamma-\gamma_{1,h}$ & $\delta_h$ & $\mathrm{ratio}$ & case & $\mathrm{eff}$
& $\dfrac{\eta_p^{\LS}\eta_d^{\LS}}{|\gamma-\gamma_{0,h}|}$
& $|\gamma-\bar\gamma_h|$ & $\eta_p^{\LS}\eta_d^{\LS}$ \\
\midrule
G--LS & P1 & $207307$ & $-2.25$e-7 & $1.00$e-7 & $-3.26$e-7 & $0.44$ & (a) & $1.44$ & $18.3$ & $6.27$e-8 & $4.12$e-6 \\
    & P2 & $202633$ & $-3.49$e-7 & $4.02$e-7 & $-7.52$e-7 & $0.87$ & (a) & $1.87$ & $306$ & $2.65$e-8 & $1.07$e-4 \\
    & P3 & $240404$ & $2.75$e-8 & $-4.93$e-8 & $7.68$e-8 & $0.56$ & (a) & $1.56$ & $47.6$ & $1.09$e-8 & $1.31$e-6 \\
    & P4 & $251245$ & $-7.47$e-7 & $5.95$e-8 & $-8.06$e-7 & $0.08$ & (a) & $1.08$ & $5.98$ & $3.44$e-7 & $4.46$e-6 \\
    & P5 & $231768$ & $1.79$e-6 & $-1.18$e-6 & $2.98$e-6 & $0.66$ & (a) & $1.66$ & $7.28$ & $3.05$e-7 & $1.31$e-5 \\
    & P6 & $240593$ & $4.24$e-6 & $-2.17$e-6 & $6.40$e-6 & $0.51$ & (a) & $1.51$ & $3.18$ & $1.04$e-6 & $1.35$e-5 \\
\midrule
LS--LS & P1 & $435744$ & $-7.13$e-8 & $1.16$e-7 & $-1.88$e-7 & $0.61$ & (a) & $1.61$ & $32.8$ & $2.26$e-8 & $2.34$e-6 \\
    & P2 & $428273$ & $-2.29$e-7 & $5.99$e-7 & $-8.28$e-7 & $0.38$ & (a) & $1.38$ & $282$ & $1.85$e-7 & $6.46$e-5 \\
\midrule
FD & P2, FD & $263169$ & $-2.29$e-6 & $6.08$e-6 & $-8.38$e-6 & $0.38$ & (a) & $1.38$ & $255$ & $1.90$e-6 & $5.84$e-4 \\
    & P2, G  & $263169$ & $-9.21$e-7 & $3.35$e-7 & $-1.26$e-6 & $0.36$ & (a) & $1.36$ & $317$ & $2.93$e-7 & $2.92$e-4 \\
\midrule
FV & P2 & $210782$ & $-2.29$e-6 & $1.52$e-5 & $-1.75$e-5 & $0.15$ & (a) & $1.15$ & $108$ & $6.45$e-6 & $2.47$e-4 \\
\bottomrule
\end{tabular}
\caption{Pairs from different methods (G--LS), on separate meshes
(LS--LS), from finite differences with G pairs on the same grid for
comparison (FD), and from finite volumes (FV), final meshes, with
$h=2^{-9}$ for FD. $N$ is the number of vertices of the common mesh, and for
LS--LS the sum of the numbers of vertices of the primal and the adjoint
mesh.}
\label{tab:indep}
\end{table}

For the G--LS pairs the results are close to those of LS and of G pairs
(Tables~\ref{tab:main}--\ref{tab:G}). On the final meshes the case is the
same, and the product exceeds the QoI gap by the factors $12.7$, $142$,
$17.0$, $5.5$, $4.4$ and $2.1$ in P1--P6, as for LS pairs except in P2,
where the factor is $123$. Along the
sequences, case (a) occurs on every mesh of P1, P3, P5 and P6, on $30$
of the $31$ meshes of P4, and in P2 from mesh $11$ on. The first ten
meshes of P2 are in case (b).

For the LS--LS pairs, the two sequences of P1 have the same numbers of vertices on every
mesh, since the problem is symmetric under $(x,y)\mapsto(-x,-y)$, which
exchanges $T_p$ and $T_d$. Case (a) occurs on all $33$ meshes, with
$\mathrm{eff}\in[1.23,1.97]$, and the product exceeds the QoI gap by factors
in $[11.5,15.2]$, as for LS pairs on a common mesh, although we compute
each pair on a mesh with $217872$ vertices and the overlay has $390757$. In P2 case (a) occurs on $30$ of
the $34$ meshes, and case (b) on meshes $12$--$15$, where $\delta_h$
changes sign, as on a common mesh (Section~\ref{sec:num_main}). Against the
total number of degrees of freedom of the two sequences, the fitted rates
over the last six meshes are $-0.98$ and $-1.07$ for $|\delta_h|$ and
$-1.00$ for $\eta_p^{\LS}\eta_d^{\LS}$, in P1 and P2.

The finite difference potentials satisfy no Galerkin relation. The
scheme takes the right-hand sides at the grid points, lumps the reaction
term and upwinds the convection term, and only its diffusion part is, up
to the factor $h^2$, the stiffness matrix of $S_1$ on this triangulation.
Case (a) occurs on all seven grids for both kinds of pairs, with
$\mathrm{eff}\in[1.38,1.59]$ for the finite difference pairs and
$\mathrm{eff}\in[1.24,1.36]$ for the G pairs. The fitted rates of
$|\gamma-\gamma_{0,h}|$, $|\gamma-\gamma_{1,h}|$ and $|\delta_h|$ against $N$ over the last four grids lie
between $-1.01$ and $-1.00$ for both, although the nodal errors of the
upwind scheme decrease only as $h$. On the final grid, $|\gamma-\gamma_{0,h}|$ is $2.5$
times and $|\gamma-\gamma_{1,h}|$ $18$ times that of the G pairs, and the product exceeds
the QoI gap by the factor $70$, against $233$. For the finite difference
pairs $\gamma_{0,h}$ is the more accurate on every grid, with
$|\gamma-\gamma_{1,h}|/|\gamma-\gamma_{0,h}|=2.66$ on the final grid, and
for the G pairs $\gamma_{1,h}$, with $|\gamma-\gamma_{0,h}|/|\gamma-\gamma_{1,h}|=2.75$. With central differences instead of
upwinding, the errors, the QoI gap and the estimators agree with those of the
G pairs to within $1\%$ from $h=2^{-6}$ on, so we do not list them.

For the FV pairs, the box method differs from the Galerkin method only where the
problem has a convection or reaction term, or data that are not constant
on the elements or on the boundary edges. In P1 and P7 its linear systems are those of the
Galerkin method, since $\bb=\bzero$, $c=0$, $f_1=g_1=0$, $\G_N=\emptyset$
and $\bff_2,\bg_2$ are constant on every element, and the results are
those of the G pairs with least-squares recovery, in P1 those of
Table~\ref{tab:G}. In
P3--P6, where it differs only in the integration of the data over the
boxes and their boundary edges, the errors, the QoI gap and the product on the final meshes agree
with those of Table~\ref{tab:G} to within $2\%$. We list only P2. There
case (a) occurs on $26$ of the $28$ meshes, all but the first two, with
$\mathrm{eff}\in[1.02,1.67]$. The fitted rates of $|\gamma-\gamma_{0,h}|$, $|\gamma-\gamma_{1,h}|$ and $|\delta_h|$ over the last six
meshes are $-1.05$, and the product exceeds the QoI gap by factors between
$11$ and $14$ on the last ten meshes, against $124$ on the final mesh of
the G pairs. As for the finite difference
pairs, $\gamma_{0,h}$ is the more accurate, on $25$ of the $28$ meshes and
by the factor $6.6$ on the final one.

\subsection{The interval and its midpoint}\label{sec:which}

Which of the two evaluators is the more accurate depends on the problem,
the method and the mesh, and the computed quantities do not tell. In P2,
$\gamma_{0,h}$ is the more accurate on $20$ of the $33$ meshes of the LS
sequence and $\gamma_{1,h}$ on the other $13$, and on the last ten meshes
$|\gamma-\gamma_{0,h}|/|\gamma-\gamma_{1,h}|$ lies between $0.49$ and
$1.14$. On the uniform grids of Table~\ref{tab:indep}, $\gamma_{0,h}$ is
the more accurate on every grid for the finite difference potentials, and
$\gamma_{1,h}$ for the Galerkin potentials. For this reason we report the
interval between $\gamma_{0,h}$ and $\gamma_{1,h}$, not one of them.

In case (a), $\gamma$ lies in the interval, and the midpoint
$\bar\gamma_h$ has the error
$$
|\gamma-\bar\gamma_h|=\tfrac12\big||\gamma-\gamma_{0,h}|-|\gamma-\gamma_{1,h}|\big|
=\tfrac12(1-\mathrm{ratio})\max(|\gamma-\gamma_{0,h}|,|\gamma-\gamma_{1,h}|).
$$
This error is at most $\frac12|\delta_h|$, and it is smaller than both
endpoint errors if and only if $\mathrm{ratio}>\frac13$. On the final
meshes of Table~\ref{tab:main}, the midpoint is more accurate than both
evaluators in P1, P2, P3, P5 and P6. In P4, where $\mathrm{ratio}=0.08$,
it is less accurate than $\gamma_{1,h}$. Along the sequences it is more
accurate than both evaluators on $30$ of the $31$ meshes of P1, $27$ of
the $33$ of P2, $30$ of the $31$ of P3, every mesh of P5 and P6, and $5$
of the $31$ meshes of P4. The bound $\frac12|\delta_h|$ exceeds the
midpoint error by the factor $(1+\mathrm{ratio})/(1-\mathrm{ratio})$,
which is $1.2$ in P4 and between $2.5$ and $6.1$ in the other problems of
Table~\ref{tab:main}. Since neither the case nor $\mathrm{ratio}$ is
computable, we report the midpoint together with the interval.

\subsection{Comparison with a dual-weighted residual computation}
\label{sec:num_dwr}

For P1, P2 and P7 we also compute the $S_2$ Galerkin approximation
$\tilde z_h$ of the adjoint problem, an approximation in an enriched space
(Section~\ref{sec:dwr}), on the meshes of Table~\ref{tab:G} and of
Section~\ref{sec:num_kellogg}, and the dual-weighted residual estimate
$\eta_{\rm DWR}=\ell(\tilde z_h-z_h)-a(u_h,\tilde z_h-z_h)
=\gamma_0(u_h,\tilde z_h)-\gamma_{0,h}$ of \eqref{etadwr}. For
$\gamma_{1,h}$ we use the G pairs in P1 and P2 and the mixed fluxes of
Section~\ref{sec:num_kellogg} in P7. In P1 and P2 the
Dirichlet data are zero or of size at most $e^{-20}$, and in P7 the
boundary term of Section~\ref{sec:setup} is below $0.3\%$ of
$|\gamma-\gamma_{0,h}|$ on the last ten meshes. Table~\ref{tab:dwr} shows the results on the final meshes.

\begin{table}[htbp]
\small
\setlength{\tabcolsep}{4.5pt}
\begin{tabular}{lrcccccc}
\toprule
& $N$ & $\gamma-\gamma_{0,h}$ & $\gamma-\gamma_0(u_h,\tilde z_h)$ & $\gamma-\gamma_{1,h}$ & $|\gamma-\bar\gamma_h|$
& $\eta_{\rm DWR}/(\gamma-\gamma_{0,h})$ & $\delta_h/(\gamma-\gamma_{0,h})$ \\
\midrule
P1 & $274402$ & $-1.31$e-7 & $5.97$e-14 & $1.24$e-7 & $3.46$e-9 & $1.00$ & $1.95$ \\
P2 & $234550$ & $-3.14$e-7 & $7.61$e-12 & $4.37$e-7 & $6.13$e-8 & $1.00$ & $2.39$ \\
P7 & $22856$  & $7.34$e-4 & $5.74$e-4 & $-6.17$e-4 & $5.87$e-5 & $0.22$ & $1.84$ \\
\bottomrule
\end{tabular}
\caption{The dual-weighted residual estimate with the $S_2$ adjoint
approximation $\tilde z_h$ and the QoI gap, for G potentials on the final
meshes of Tables~\ref{tab:G} and~\ref{tab:kellogg}.}
\label{tab:dwr}
\end{table}

In P1 and P2, $\eta_{\rm DWR}/(\gamma-\gamma_{0,h})=1.00$ on the last ten
meshes, because $\gamma_0(u_h,\tilde z_h)$ is much more accurate than
$\gamma_{0,h}$. In P2 it is more accurate by the factor $4\cdot10^4$, and
its error decreases with the fitted rate $-1.95$, against $-1.03$ for
$|\gamma-\gamma_{0,h}|$. In P1 it reaches the accuracy of the reference
value. The QoI gap is in case (a) in both problems and estimates
the larger of the errors of $\gamma_{0,h}$ and $\gamma_{1,h}$. In P7,
$\eta_{\rm DWR}/(\gamma-\gamma_{0,h})$ lies in $[0.22,0.23]$ on the last
ten meshes, so that the error of $\gamma_0(u_h,\tilde z_h)$ is $78\%$ of
that of $\gamma_{0,h}$. The adjoint solution has in general, like $u$, the
singular behavior $r^{0.1}$ at the origin, and on a given mesh a higher
polynomial degree does little to improve the approximation of such a
singularity. The QoI gap is in case (a) on every mesh of P7.

\subsection{Dominating reaction and the reaction-restored evaluator \texorpdfstring{$\gamma_2$}{gamma2}}\label{sec:num_sweep}

We test the reaction-restored evaluator $\gamma_2$ of \eqref{gamma2} and
the reaction-restored QoI gap on P2 with $\epsilon=0.05$. We vary the
reaction coefficient from $c=1$ to $c=10^8$ and keep the exact solutions
$u$ and $z$ fixed, so that $\gamma=\gamma|_{c=0}+c\,(z,u)$. We compute the least-squares
approximations with the balance terms weighted by $c^{-1/2}$, as in
Remark~\ref{rem:weight}, and refine the meshes with the balanced
indicator of the weighted estimators. For each $c$ we take all quantities
from one mesh of its adaptive sequence, the one whose number of
vertices is closest to $3140$, the smallest final number of vertices
among the nine sequences of $28$ meshes. As $c$ grows, the parameter
$\vartheta=|\bb|/\sqrt{\epsilon c}$ of Proposition~\ref{prop:delta2}
decreases from $\vartheta\gg1$ to $\vartheta\ll1$, where reaction
dominates.

The data $f_1$ and $g_1$ grow linearly in $c$, and so do the integrals of
data in the evaluators, for instance $(c^{-1}f_1,g_1)$ in $\gamma_2$, which
is of size $c\,(u,z)$. Their quadrature errors also grow linearly in
$c$. Let $Q(uz)$ be the integral of $uz$ computed by the three-point rule
of degree two. At $c=10^8$, on the mesh of Table~\ref{tab:sweep},
$c\,|(u,z)-Q(uz)|$ is $0.51$, more than three orders of magnitude larger
than the error of $\gamma_{2,h}$. For this reason we compute all
evaluators with the $19$-point rule of degree nine. Replacing it by the
rule of degree eight changes them by less than
$10^{-3}|\gamma-\gamma_{2,h}|$.

\begin{table}[htbp]
\footnotesize
\setlength{\tabcolsep}{3pt}
\begin{tabular}{ccrccccccc}
\toprule
$c$ & $\vartheta$ & $N$ & $|\gamma-\gamma_{0,h}|$ & $|\gamma-\gamma_{1,h}|$
& $|\gamma-\gamma_{2,h}|$ & $|\delta_h|$
& $|\delta_{2,h}|$ & $\dfrac{|\delta_{2,h}|}{|\gamma-\gamma_{0,h}|}$ & $\eta_{p,c}^{\LS}\eta_{d,c}^{\LS}$ \\
\midrule
$1$ & $16.1$ & $3348$ & $3.42$e-5 & $2.54$e-5 & $1.79$e-3 & $8.77$e-6 & $1.76$e-3 & $51.3$ & $1.55$e-2 \\
$10$ & $5.10$ & $3509$ & $1.27$e-6 & $1.08$e-4 & $3.24$e-4 & $1.10$e-4 & $3.25$e-4 & $255$ & $1.20$e-2 \\
$10^{2}$ & $1.61$ & $2797$ & $1.53$e-5 & $1.62$e-4 & $2.00$e-4 & $1.78$e-4 & $2.15$e-4 & $14.0$ & $1.50$e-2 \\
$10^{3}$ & $0.510$ & $3183$ & $1.71$e-5 & $1.61$e-4 & $1.39$e-4 & $1.44$e-4 & $1.22$e-4 & $7.14$ & $1.32$e-2 \\
$10^{4}$ & $0.161$ & $3425$ & $4.90$e-5 & $1.87$e-4 & $1.28$e-4 & $1.38$e-4 & $7.90$e-5 & $1.61$ & $1.46$e-2 \\
$10^{5}$ & $0.0510$ & $2761$ & $9.49$e-4 & $1.17$e-3 & $2.11$e-4 & $2.19$e-4 & $7.38$e-4 & $0.778$ & $5.04$e-2 \\
$10^{6}$ & $0.0161$ & $3197$ & $8.43$e-3 & $8.65$e-3 & $2.14$e-4 & $2.17$e-4 & $8.22$e-3 & $0.975$ & $2.61$e-1 \\
$10^{7}$ & $5.10\cdot10^{-3}$ & $3140$ & $8.31$e-2 & $8.33$e-2 & $2.12$e-4 & $2.16$e-4 & $8.29$e-2 & $0.997$ & $2.58$ \\
$10^{8}$ & $1.61\cdot10^{-3}$ & $3143$ & $8.34$e-1 & $8.34$e-1 & $2.20$e-4 & $2.24$e-4 & $8.34$e-1 & $1.00$ & $25.6$ \\
\bottomrule
\end{tabular}
\caption{P2 with $\epsilon=0.05$ and $c$ from $1$ to $10^8$, least
squares with the balance terms weighted by $c^{-1/2}$. For each $c$ we
use one mesh, with $N$ vertices, and
$\delta_{2,h}=\gamma_{2,h}-\gamma_{0,h}$.}
\label{tab:sweep}
\end{table}

Table~\ref{tab:sweep} shows two ranges of $c$ with different behavior. For $c\le10^2$, where
$\vartheta\ge1.6$, $\gamma_{2,h}$ has the largest error of the three
evaluators, by the factor $70$ against $\gamma_{1,h}$ at $c=1$. There
$|\delta_{2,h}|$ is close to the error of $\gamma_{2,h}$ and exceeds
$|\gamma-\gamma_{0,h}|$ by factors between $14$ and $255$.
For $c\ge10^5$ the errors of
$\gamma_{0,h}$ and $\gamma_{1,h}$ grow in proportion to $c$, from
$9.5\cdot10^{-4}$ at $c=10^5$ to $0.83$ at $c=10^8$, and for $c\ge10^6$
they nearly coincide. The QoI gap $|\delta_h|$ stays at
$2.2\cdot10^{-4}$, because the signed QoI gap does not contain the
reaction term $(ce_p,e_d)$ (Section~\ref{sec:gamma2}). With
$\mathrm{TOL}=10^{-3}$, the criterion $|\delta_h|\le\mathrm{TOL}$ holds for
every $c\ge10^5$, although the errors of $\gamma_{0,h}$ and $\gamma_{1,h}$
reach $0.83$. For this reason we use $|\delta_{2,h}|$ where reaction
dominates (Section~\ref{sec:procedure}).
The error of
$\gamma_{2,h}$ does not grow. It lies between $1.28\cdot10^{-4}$ and
$2.20\cdot10^{-4}$ for all $c\ge10^3$, and at $c=10^8$ it is smaller than
that of $\gamma_{0,h}$ by the factor $3.8\cdot10^3$.

For large $c$, the reaction-restored QoI gap therefore estimates the
error of $\gamma_{0,h}$. For $c=10^6,10^7,10^8$ the quotient
$|\delta_{2,h}|/|\gamma-\gamma_{0,h}|$ is $0.975$, $0.997$ and $0.9997$,
and $\mathrm{ratio}=|\gamma-\gamma_{2,h}|/
|\gamma-\gamma_{0,h}|$ is $0.025$, $0.0026$ and $0.00026$, so that
$(1-\mathrm{ratio})|\gamma-\gamma_{0,h}|\le|\delta_{2,h}|
\le(1+\mathrm{ratio})|\gamma-\gamma_{0,h}|$ by
Proposition~\ref{prop:twosided}. For $c\ge10^6$, where
$\vartheta\le0.016$, $|\delta_{2,h}|$ is $0.03$ times
$\eta_{p,c}^{\LS}\eta_{d,c}^{\LS}$, so the bound \eqref{delta2bound} is
pessimistic by a factor of about $30$. Between these two ranges, at
$c=10^3$ and $10^4$, $\gamma_{0,h}$ is the most accurate, by the factors
$8.1$ and $2.6$ over $\gamma_{2,h}$, and $|\delta_{2,h}|$ exceeds
$|\gamma-\gamma_{0,h}|$ by the factors $7.1$ and $1.6$.

Over all values of $c$, $|\delta_{2,h}|$ lies between $0.62$ and $1.08$
times $\max(|\gamma-\gamma_{0,h}|,|\gamma-\gamma_{2,h}|)$. The pair
$\gamma_{0,h},\gamma_{2,h}$ is in case (a) at $c=10$ and $10^2$ and in
case (b) at the other values of $c$, with $\mathrm{ratio}\le0.38$.

\section{Concluding remarks}\label{sec:conc}

The product of the primal and adjoint least-squares estimators is a
reliable bound for the error in the quantity of interest, and its local
contributions are good refinement indicators. It is, however, too
pessimistic for deciding when to stop. In our experiments it exceeds the
larger of the errors of the two evaluators by factors of up to several
hundred, and if we stopped with it we would refine far more than
needed.

We compute the signed QoI gap from the same residuals. Like the dual-weighted residual estimate, it is the
difference of two values of the quantity of interest, but the second
value comes from the fluxes and not from a more accurate approximation.
In our
experiments the QoI gap is of the size of the larger of the errors of the
two evaluators. We therefore refine
with the least-squares estimators and stop with the QoI gap. Where reaction
dominates, we stop with the reaction-restored QoI gap instead.

With the least-squares flux recovery, we can use the QoI gap with conforming
potentials from any method. Our experiments include Galerkin, finite
difference and finite volume potentials. We leave the boundary term that
appears when we interpolate the Dirichlet data to the next paper of this
series.

\end{document}